\documentclass[12pt]{article}

\usepackage{indentfirst}

\usepackage[latin1]{inputenc}
\usepackage{amsmath,amssymb}
\usepackage{amsmath,accents}
\usepackage{alltt,graphics,graphpap,color}
\usepackage[dvips]{graphicx}
\usepackage{amsfonts}
\usepackage{amscd}
\usepackage{mathrsfs}
\usepackage{hyperref}

\hypersetup{linktoc=page,colorlinks=true,linkcolor=red,  citecolor=blue, }

\newtheorem{Theorem}{Theorem}[section]
\newtheorem{Thdefi}[Theorem]{Theorem / Definition}

\newtheorem{Proposition}[Theorem]{Proposition}
\newtheorem{Lemma}[Theorem]{Lemma}
\newtheorem{Corollary}[Theorem]{Corollary}
\newtheorem{Remark}[Theorem]{Remark}
\newtheorem{Example}[Theorem]{Example}
\newtheorem{Notation}[Theorem]{Notation}

\newcommand{\ch}{\mathbb {CH} }
\newcommand{\sss}{\mathbb S}
\newcommand{\coq}{\cosh^2(\rho)}

\newcommand{\siq}{\sinh^2(\rho)}
\newcommand{\sih}{\sinh(\rho)}
\newcommand{\coh}{\cosh(\rho)}
\newcommand{\coqt}{\cosh^2(\rho(t))}
\newcommand{\siqt}{\sinh^2(\rho(t))}

\newcommand{\mm}{\mathcal {M}}

\newcommand{\parho}{\frac{\partial}{\partial\rho}}
\newcommand{\pat}{\frac{\partial}{\partial t}}

\newcommand{\proof}{\noindent\emph{Proof. }}
\newcommand{\cvd}{\hfill$\square$ \bigskip}

\newcommand{\ep}{\mu}

\newcommand{\ns}{{\nabla_{\sigma}}}
\newcommand{\nid}{\noindent}
\begin{document}

%\def\qed{\hbox{\hskip 6pt\vrule width6pt height7pt
%depth1pt  \hskip1pt}\bigskip}

%\setlength\parindent{0pt}
%\vspace*{1 cm}

\title{Inverse mean curvature flow in complex hyperbolic space }

\author{\sc Giuseppe Pipoli\footnote{This research is supported by the ERC Avanced Grant 320939, Geometry
and Topology of Open Manifolds (GETOM)}}
\date{}

\maketitle

%\vspace{3cm}
%\textbf{\Large 6th version}
\vspace{1,5cm}

{\small \noindent {\bf Abstract:} We consider the evolution by inverse mean curvature flow of a closed, mean convex and star-shaped hypersurface in the complex hyperbolic space. We prove that the flow is defined for any positive time, the evolving hypersurface stays star-shaped and mean convex. Moreover the induced metric converges, after rescaling, to a conformal multiple of the standard sub-Riemannian metric on the sphere. Finally we show that there exists a family of examples such that the Webster curvature of this sub-Riemannian limit is not constant.} 

%\medskip

%{\small \noindent {\bf R\' esum\' e :} Nous consid\' erons l'\' evolution par l'inverse de la courbure moyenne d'une surface \' etoil\' ee, ferm\' ee et \` a courbure moyenne positive dans l'espace hyperbolique complexe. Nous montrons que le flot est d\' efini pour tout temps positif et que la surface reste \' etoil\' ee et \` a courbure moyenne positive. De plus, la m\' etrique induite, apr\` es un changement d'\' echelle, converge vers un multiple conforme de la m\' etrique sous-riemannienne standard sur la sph\` ere de dimension impaire. Nous allons montrer l'existence d'exemples de donn\' ees initiales telles que cette limite sous-riemannienne n'a pas courbure de Webster constante.}

\medskip


%\noindent {\bf Mots-clefs :} Flot par l'inverse de la courbure moyenne, espace hyperbolique complexe, g\' eom\' etrie sous-Riemannienne.\\

 \bigskip

\section{Introduction} \setcounter{equation}{0}

During last years geometric flows of submanifolds of Riemannian manifolds have been studied intensively. In the class of expanding flows, the leading example is the inverse mean curvature flow. In this paper we consider the evolution by inverse mean curvature flow of real hypersurfaces of the complex hyperbolic space $\ch^n$, with $n\geq 2$. For any given smooth hypersurface $F_0:\mm\rightarrow\ch^n$ , the solution of the inverse mean curvature flow with initial datum $F_0$ is a one-parameter family of smooth immersion $F:\mm\times\left[0,T\right)\rightarrow\ch^n$ such that
\begin{equation}\label{imcf_def}
\left\{\begin{array}{rcl}
\displaystyle{\frac{\partial F}{\partial t}} & = & \displaystyle{\frac 1H \nu,}\\
F(\cdot,0) & = & F_0,
\end{array}\right.
\end{equation}
where $H$ is the mean curvature of $F_t=F(\cdot,t)$ and $\nu$ is the outward unit normal vector of $\mm_t=F_t(\mm)$. It is the main tool used in the celebrated paper of Huisken and Ilmanen in \cite{HI} for proving the Penrose inequality.

In this paper we restrict to the class of star-shaped hypersurface in $\ch^n$. The inverse mean curvature flow of star-shaped hypersurfaces has been already studied in different ambient manifolds: for example the Euclidean space \cite{Ge1,Ur}, the hyperbolic space \cite{Ge3,HW}, asymptotic hyperbolic spaces \cite{Ne}, rotationally symmetric spaces \cite{Di} and warped products \cite{Sc,Zh}. In any case it was proved that the flow is defined for any positive time and the evolving hypersurface stays star-shaped for all the life of the flow. Inverse mean curvature flow of star-shaped hypersurfaces in rank one symmetric spaces have been considered for different pourposes in \cite{KS} too. 

The geometry of the ambient manifold influences the nature of the limit of the induced metric. In fact Gerhardt \cite{Ge1} and and Urbas \cite{Ur} proved indipendently that, for any star-shaped hypersurface of the Euclidean space, the limit metric is, up to rescaling, always the standard round metric on the sphere. In \cite{HW} P.K. Hung and M.T. Wang showed that, when the ambient manifold is the hyperbolic space, the limit metric is not always round. More precisely it is a conformal multiple of the standard round metric on the sphere and so it is round only in special cases. The case studied in the present paper has some similarities with the previous results, but a new phenomenon appears: even after rescaling, the evolving metric blows up along a direction. Hence the limit metric is no more Riemannian, but only sub-Riemannian: it is defined only on a codimension-one distribution. The main theorem proved is the following.

\begin{Theorem}\label{main}
For any $\mm_0$ closed, mean convex and $\sss^1$-invariant star-shaped hypersurface in $\ch^n$ let $\mm_t$ be its evolution by inverse mean curvature flow, let $g_t$ be the induced metric on $\mm_t$ and $\theta_t$ the induced contact form. Then:
\begin{itemize}
\item[(1)] $\mm_t$ is $\sss^1$-invariant, star-shaped and mean convex for any time $t$;
\item[(2)] the flow is defined for any positive time;
\item[(3)] there is a smooth $\sss^1$-invariant function $f$ such that the rescaled induced metric $\tilde g_t=\left|\mm_t\right|^{-\frac 1n}g_t$ converges to a sub-Riemannian metric $\tilde g_{\infty}=e^{2f}\sigma_{sR}$ (i.e. a conformal multiple of the standard sub-Riemannian metric on the sphere $\sss^{2n-1}$) and the rescaled contact form $\tilde{\theta}=\left|\mm_t\right|^{-\frac 1n}\theta_t$ converges to $\tilde{\theta}_{\infty}=e^{2f}\hat{\theta}$, where $\hat{\theta}$ is the standard contact form on the odd dimensional sphere;
\item[(4)] moreover there are examples of $\mm_0$ such that the limit does not have constant Webster scalar curvature.

\end{itemize}
%$\mm_t$ is defined, mean convex and star-shaped for any positive time. The rescaled induced metric $\tilde g_t=\left|\mm_t\right|^{-\frac 1n}g_t$ converges as $t\rightarrow\infty$ to $\tilde g_{\infty}=e^{2f}\sigma_{sR}$ -a conformal multiple of the standard sub-Riemannian metric on the sphere $\sss^{2n-1}$ -  where $f$ is a smooth function. Moreover there are examples of $\mm_0$ such that $\tilde g_{\infty}$ does not have constant Webster curvature.
\end{Theorem}

After proving that the flow can be extended for any positive time, we have that the volume of $\mm_t$ becomes arbitrary large and then $\mm_t$ ``explores'' the structure at infinity of the ambient manifold as $t$ tends to infinity. Hence the convergence to a sub-Riemannian metric is not surprising because the boundary at infinity of $\ch^n$ is $(\sss^{2n-1}, \sigma_{sR})$, the one point compactification of the Heisenberg group of dimension $2n-1$ endowed with its standard sub-Riemannian metric. Different initial data explore this structure at infinity in different ways, but our result shows that we remain in the conformal class of $\sigma_{sR}$.

Obviously for any finite time $t$, $\tilde g_t$ is a Riemannian metric, but there is a direction in which the metric is blowing up. This special direction is $J\nu$, where $J$ is the complex structure of $\ch^n$ and $\nu$ is the outward unit normal vector field of $\mm_t$. Note that, since we are considering only submanifolds of codimension one, $J\nu$ is for sure tangent to $\mm_t$. One of the main difficulties in generalizing the previous results is to describe the contribution of this special direction.
%part $(1)$ and $(2)$ of Theorem \ref{main} 
Its presence gives also a new phenomenon not present in the previous literature. The second fundamental form converges to that of an horosphere with an exponential speed but, unlike for example \cite{Ge1,Sc,Zh}, we have that the speed is not the same for any initial datum: very symmetric hypersurfaces converges twice as fast as the generic $\sss^1$-invariant hypersurface (see Remark \ref{remark} below for more details).

If we try to study the limit of the sectional curvature of this family of metric $\tilde{g}_t$, it always diverge: this is a general behaviour when we try to approximate a sub-Riemannian metric with a family of Riemannian metrics. For that reason an other notion of curvature is required. We will use in particular the Webster curvature.

It is very easy to find hypersurfaces of $\ch^n$ such that $\tilde{g}_{\infty}$ has constant Webster curvature. It is the case of the geodesic spheres: as we will see in detail in Section \ref{sez sf geo}, the evolution of a geodesic sphere is a family of geodesic spheres and the function $f$ is constant. On the other hand, the search for an example with a non-trivial limit is much more challenging.
The main tool for studying the roundness of the limit is the following Brown-York like quantity: for any star-shaped hypersurface $\mm$
$$
Q(\mm)=|\mm|^{-1+\frac 1n}\int_{\mm}\left(H-\hat{H}\right)d\mu,
$$
where $|\mm|$ is the volume of $\mm$ and, if $\rho$ is the radial function defining $\mm$, $\hat{H}$ is the value of the mean curvature of a geodesic sphere of radius $\rho$ (see \eqref{hat_H} for the explicit definition). $Q$ gives a measure of how $\mm$ is far to being a geodesic sphere. It is not a measure in the strict sense because $Q$ has not a sign and, even if it is zero for geodesic spheres, it is not in general true the opposite. In the final section of this paper we found the desired non-trivial examples estimating the behaviour of $Q$ along the inverse mean curvature flow.

 %We found such an example in the class of $\sss^1$-invariant submanifolds, introducing a quantity $Q(\mm)$ which gives a measure of how $\mm$ is far to being a geodesic sphere (see \eqref{Q_def} for a precise definition) and estimating its behaviour along the inverse mean curvature flow.

This paper is organized as follows. In Section \ref{sez 2} we collect some preliminaries and we fix some notations. In Section \ref{sez 3} we compute the main geometric quantities for a star-shaped hypersurface in $\ch^n$, like the induced metric, the second fundamental form and the mean curvature. In Section \ref{sez sf geo} we have a simple but meaningful example, i.e. the evolution of the geodesic spheres. In Section \ref{sez 5} we estimate the norm of the gradient of the radial function. As consequence we have that the property of being star-shaped and the mean convexity are preserved by the flow. The study of the derivatives of the radial function continues in Sections \ref{sez 6} and \ref{sez 7}. In particular we prove that the solution of the flow is defined for any positive time. Section \ref{sez 8} is devoted to the prove of the convergence of the rescaled induced metric to a sub-Riemannian limit. In the last Section we conclude the proof of Theorem \ref{main} studying the Webster curvature of the limit metric and giving a family of non-trivial examples.

Finally we would like to announce that the ideas developed in the present paper have been extended in \cite{Pi2} in the case of the next rank one symmetric space, that is the quaternionic hyperbolic space. An analogous of Theorem \ref{main} holds in this other setting too.

\medskip

\noindent\textbf{Acknowledgements:} The author of the present paper would like to thank G\' erard Besson for many useful conversations had during the course of this research and his constant interest and encouragement in this project.

\section{Preliminaries}  \setcounter{equation}{0}\label{sez 2}

\subsection{Riemannian and sub-Riemannian metrics on the sphere}
Every hypersurface considered in this paper is closed and star-shaped and so it is an embedding of $\sss^{2n-1}$, the sphere of dimension $2n-1$, into $\ch^n$. On that sphere we have different ``standard'' metrics. Let $\sigma$ be the usual Riemannian metric on $\sss^{2n-1}$ with constant sectional curvature equal to 1. Since the dimension is odd, we can distinguish an important vector field: we can think the sphere embedded in $(\mathbb C^n, J)$, then, if $\nu$ is the unit outward normal to $\sss^{2n-1}$, $\xi=J\nu$ is an unit tangent vector field on the sphere. It is often called the \emph{Hopf vector field}, because it can be characterized also as the unit vector field tangent to the fibers of the Hopf fibration: $\pi:\sss^{2n-1}\longrightarrow\mathbb{CP}^{n-1}$. It allows us to define the \emph{horizontal distribution} 
\begin{equation}\label{horizontal distribution}
\mathcal{H}=\left\{X\in T(\sss^{2n-1})\left|\sigma(X,\xi)=0\right.\right\}
\end{equation} 

The \emph{Berger metrics} is a family of deformations of $\sigma$ in the direction of $\xi$: for any $\mu>0$ let $e_{\mu}$ be the Riemannian metric defined by
$$
\left\{\begin{array}{rccl}
e_{\mu}(X,Y) & = & \sigma(X,Y) &\quad\text{for any } X,\ Y\in \mathcal H;\\
e_{\mu}(X,\xi) & = & 0 &\quad\text{for any } X\in \mathcal H;\\
e_{\mu}(\xi,\xi) & = & \mu.
\end{array}\right.
$$
When $\mu$ converges to infinity, the metric $e_{\mu}$ degenerates on the directions tangent to $\xi$. At the limit we get $\sigma_{sR}$ the \emph{standard sub-Riemannian metric}, it is defined only on $\mathcal H$, but, since $\mathcal H +\left[\mathcal H,\mathcal H\right]$ is the whole tangent space, $\sigma_{sR}$ is enough to define a distance on $\sss^{2n-1}$, called the \emph{Carnot-Caratheodory distance}. The following comparison between the Levi-Civita connection of $e_{\ep}$ and that of $\sigma=e_1$ will be used below.

\begin{Lemma}\label{diff_LCC}
Fix a $\sigma$-orthonormal basis $(Y_1,\cdots,Y_{2n-1})$ of $\sss^{2n-1}$ such that $Y_1=\xi$ and for every $r$ $Y_{2r+1}=JY_{2r}$. Let us denote with ${\nabla}_e$ ($\nabla_{\sigma}$ respectively) the Levi-Civita connection associated to the metric $e_{\mu}$ ($\sigma$ respectively). Then for every $1\leq i,\ j\leq 2n-1$ we have:
$$
\nabla_{e\ Y_i}Y_j-\nabla_{\sigma\ Y_i}Y_j = \left\{\begin{array}{ll}
(1-\ep)JY_j & \text{if }i=1\neq j; \\
 (1-\ep)JY_i & \text{if }j=1\neq i;\\
0 & \text{otherwise.}
\end{array}\right.
$$
\end{Lemma}
\proof Obviously $\nabla_{e\ Y_1}Y_1=\nabla_{\sigma\ Y_1}Y_1=0$ since the Hopf vector field is tangent to the fibers of the Hopf fibration, and they are geodesic for every $\ep$. The metric $e_{\ep}$ can be seen as the metric on the total space of the canonical variation of parameter $\ep$ of the Hopf fibration. By Lemma 3 of \cite{O} and Lemma 9.69 in \cite{Be} we have:
$$
\begin{array}{ll}
\nabla_{e\ Y_i}Y_j=\nabla_{\sigma\ Y_i}Y_j & \text{ if } i,j\neq 1;\\
\nabla_{e\ Y_i}Y_1= \ep \nabla_{\sigma\ Y_i}Y_1=-\ep JY_i & \text{ if } i\neq 1,
\end{array}
$$
where the last equality comes from an explicit computation. Finally, if $i\neq 1$, we get:
$$
\nabla_{e\ Y_1}Y_i=\nabla_{e\ Y_i}Y_1+\left[Y_1,Y_i\right] = (1-\ep)JY_i+\nabla_{\sigma\ Y_1}Y_i.
$$
\cvd

\begin{Notation}\label{notazioni hess}
We introduce the following notation in order to distinguish between derivatives of a function with respect to different metrics. For any given function $f:\sss^{2n-1}\rightarrow\mathbb{R}$, let $f_{ij}$ ($\hat{f}_{ij}$ respectively) be the components of the Hessian of $f$ with respect to $\sigma$ ($e_{\ep}$ respectively). The value of $\ep$ will be clear from the context. The indices go up and down with the associated metric: for instance $\hat{f}_i^k=\hat{f}_{ij}e_{\ep}^{jk}$, while $f_i^k=f_{ij}\sigma^{jk}$. Analogous notations will be used for higher order derivatives.
\end{Notation}

On the sphere of odd dimension $\sss^1$ acts in the following way:
\begin{equation}\label{azione s1}
\begin{array}{llll}
S: & \sss^1\times\sss^{2n-1}\subset\mathbb{C}\times\mathbb{C}^n& \rightarrow & \sss^{2n-1}\subset\mathbb{C}^n\\
& (e^{i\vartheta},z_1,\cdots,z_n) & \mapsto & (e^{i\vartheta}z_1,\cdots,e^{i\vartheta}z_n).
\end{array}
\end{equation}
The action is by isometries for any Berger metric. A function $\varphi:\sss^{2n-1}\rightarrow\mathbb{R}$ is said $\sss^1$-invariant if it is invariant under the action of $S$.

\begin{Lemma}\label{hessiani}
Let $\varphi$ be an $\sss^1$-invariant smooth function. With respect to the basis introduced in the previous Lemma, the Hessian of $\varphi$ with respect to $e_{\ep}$ is:
$$
\hat{\varphi}_{ij}=\left(\begin{array}{cc}
0 & \ep JY_j(\varphi)\\
\ep JY_i(\varphi)& {\varphi}_{ij}
\end{array}\right),
$$
where we are using Notations \ref{notazioni hess} and $J$ is the complex structure of $\mathbb{C}^n$ (once again we are considering $\sss^{2n-1}$ embedded in $\mathbb{C}^n$). Taking the trace and the norm of the Hessian, in particular we have:
\begin{eqnarray*}
\Delta_e\varphi&=&\Delta_{\sigma}\varphi;\\
|\nabla^2_{e}\varphi|^2_{e} &=&|\nabla^2_{\sigma}\varphi|^2_{\sigma}+2(\ep-1)|\nabla_{\sigma}\varphi|^2_{\sigma}.
\end{eqnarray*}
\end{Lemma}
\proof
Since $\varphi$ is $\sss^1$-invariant we have that $\varphi_1=Y_1(\varphi)=0$.  From the previous Lemma we get:
\begin{eqnarray*}
\hat{\varphi}_{11}& =& Y_1Y_1(\varphi)-\nabla_{e\ Y_1}Y_1(\varphi)=0.
%(\nabla^2_{e}\varphi)_{1j}& =& Y_jY_1(\varphi)-\nabla_{e\ Y_j}Y_1\varphi=\ep JY_i(\varphi)\\
\end{eqnarray*}
For every $i\neq 1$ we can compute:
\begin{eqnarray*}
%(\nabla^2_{e}\varphi)_{11}& =& Y_1Y_1(\varphi)-\nabla_{e\ Y_1}Y_1\varphi=0;\\
\hat{\varphi}_{1i}& =& Y_iY_1(\varphi)-\nabla_{e\ Y_i}Y_1(\varphi)=\ep JY_i(\varphi);\\
\hat{\varphi}_{i1}& =& Y_1Y_i(\varphi)-\nabla_{e\ Y_1}Y_i(\varphi)\\
&=& \left[ Y_i,Y_1\right](\varphi) +Y_iY_1(\varphi)-\nabla_{e\ Y_1}Y_i(\varphi)\\
&=&\hat{\varphi}_{1i}.
\end{eqnarray*}
Moreover, if both index are different from 1, we have:
\begin{eqnarray*}
\hat{\varphi}_{ij}& =& Y_jY_i(\varphi)-\nabla_{e\ Y_j}Y_i(\varphi)\\
& = & Y_jY_i(\varphi)-\nabla_{\sigma\ Y_j}Y_i(\varphi)={\varphi}_{ij}.
\end{eqnarray*}
We point out that, as a consequence of the symmetries considered, 
$$
|\nabla_e\varphi|^2_e=|\nabla_{\sigma}\varphi|^2_{\sigma}
$$
Taking into account this remark, the formulas for the Laplacian and the norm of the Hessian follow after some trivial computations.
\cvd

\subsection{Webster curvature}

The Webster curvature is a notion of CR geometry. In this section we consider only the case of the sphere, what follows can be said with much more generality. We refer to the monograph \cite{DT} for all the details. Fix the standard CR structure on $\sss^{2n-1}$, then $\mathcal{H}$, defined in \eqref{horizontal distribution}, is the horizontal distribution of this structure. On $\mathcal{H}$ we have a complex structure $J$. Let $\pi:T(\sss^{2n-1})\rightarrow\mathcal{H}$ be the canonical projection. For any $1$-form $\theta$ such that $\ker\theta=\mathcal{H}$, there exists a vector field $\xi_{\theta}$ such that $\theta(\xi_{\theta})=1$ and $d\theta(\xi_{\theta},\cdot)=0$. We can extend $J$ to the whole tangent space of $\sss^{2n-1}$ requiring that $J\xi_{\theta}=0$. For every $X,Y\in\mathcal{H}$ we can define $G_{\theta}(X,Y)=d\theta(X,JY)$. The metric $g_{\theta}=\pi G_{\theta}+\theta^2$ is called Webster metric associated to $\theta$. 

\begin{Thdefi}
Let $\theta$, $J$, $\xi_{\theta}$ and $g_{\theta}$ be as before. There exists an unique linear connection $\nabla^{TW}$ such that $\nabla^{TW} J=\nabla^{TW}\theta=\nabla^{TW}\xi_{\theta}=\nabla^{TW}g_{\theta}=0$ and with torsion $T$ of pure type, i.e. for every $X,Y\in\mathcal{H}$ the torsion satisfies:
$$
\left\{\begin{array}{l}
T(X,Y)=d\theta(X,Y)\xi_{\theta};\\
T(X,\xi_{\theta})\in\mathcal{H}\\ 
g_{\theta}(T(X,\xi_{\theta}),Y)=g_{\theta}(T(Y,\xi_{\theta}),X)=-g(T(JX,\xi_{\theta}),JY).
\end{array}\right.
$$
This connection $\nabla^{TW}$ is called the Tanaka-Webster connection associated to $\theta$.
\end{Thdefi}

\nid The Webster curvature of $\theta$ is the curvature defined in the usual way, but using the Tanaka-Webster connection instead of the Levi-Civita connection. 

On $\sss^{2n-1}$ we have a standard contact form: $\hat{\theta}(\cdot)=\sigma(\xi,\cdot)$. Whit respect to this form we have: $\xi_{\hat{\theta}}=\xi$, $G_{\hat{\theta}}=\sigma_{sR}$ and $g_{\hat{\theta}}=\sigma$. Obviously the Webster curvature of $\hat{\theta}$ is constant. A metric of the form $e^{2f}\sigma_{sR}$ can be thought as the restriction to $\mathcal{H}\times\mathcal{H}$ of the Webster metric associated to the $1$-form $e^{2f}\hat{\theta}$. Then we will talk indifferently about the  Webster curvature of a sub-Riemannian metric or of a contact form. 

Once we defined the appropriate notion of curvature in CR geometry, a natural question is the Yamabe problem for CR manifolds. This problem was solved in a great generality, but in our special case it can be reformulated as follows: what are the functions $f$ such that the conformal multiple $e^{2f}\sigma_{sR}$ has constant Webster scalar curvature? The answer is given by the following formula by Jerison and Lee \cite{JL}: the Webster scalar curvature of $e^{2f}\sigma_{sR}$ is constant if and only if there are $u>0$, $c>0$ and $\zeta\in\sss^{2n-1}$ such
\begin{equation}\label{JL}
e^{-2f}(z)=c\left|\cosh(u)+\sinh(u)z\cdot\bar{\zeta}\right|^2, \quad\forall z\in\sss^{2n-1}.
\end{equation}
Here we are considering the odd-dimensional sphere immersed in $\mathbb C^n$ and the norm and the product are the usual ones in $\mathbb C^n$.

Since in this paper we are considering $\sss^1$-invariant hypersurface, we show that, in presence of this kind of symmetry, the Jerison and Lee's formula becomes much simpler.
%The formula of Jerison and Lee \eqref{JL} suggests that, in presence of symmetries, the problem can be simplified. 

\begin{Lemma}\label{caratterizzazione_curv}
Let $f:\sss^{2n-1}\rightarrow\mathbb R$ be an $\sss^1$-invariant function. The following are equivalent:
\begin{itemize}
\item[(a)] $f$ satisfies \eqref{JL},
%\item[(b)] $e^{-f}$ is a linear combination of constants and first eigenfunctions of $\sss^{2n-1}$,
\item[(b)] $f$ is constant.
\end{itemize}
\end{Lemma}
\proof If $f$ is constant, $(a)$ holds trivially.  In order to prove the opposite implication, let  $u>0$ and $\zeta=(a_j-ib_j)_{j=1,\dots,n}\in\sss^{2n-1}\subset\mathbb C^n$ be as in equation \eqref{JL}. By hypothesis we have that for any $z\in\sss^{2n-1}$ and for any $\vartheta=x+iy\in\sss^1\subset\mathbb C$ we have

\begin{equation}\label{eq080}
\left|\cosh(u)+\sinh(u)(e^{i\vartheta}z)\cdot\bar{\zeta}\right|^2=\left|\cosh(u)+\sinh(u)z\cdot\bar{\zeta}\right|^2
\end{equation}

\noindent Fix $z_1=(1,0,\dots,0)$, then $e^{i\vartheta}z\cdot\bar{\zeta}=xa_1-yb_1+i(xb_1+ya_1)$, while $z\cdot\bar{\zeta}=a_1+ib_1$. Specifying \eqref{eq080} for $z_1$ we have, after some trivial computations, that 
$$
\sinh(u)\cosh(u)a_1=\sinh(u)\cosh(u)(xa_1-yb_1),\quad \forall x,y \text{ such that }x^2+y^2=1
$$
then $\sinh(u)=0$ or $a_1=b_1=0$. For any $k=2,\dots,n$ we can repeat the same computations for $z_k=(\delta_{kj})_{j=1,\dots,n}$ and we get that $\sinh(u)=0$ or $a_k=b_k=0$. Since $\zeta\neq0$ the only possibility is that $\sinh(u)=0$ and then we have the thesis.
\cvd

\subsection{Complex hyperbolic space}
The complex hyperbolic space is the complex analogous of the real hyperbolic space. It can be defined in many equivalent ways, but, for our purpose, it is convenient to introduce polar coordinates. Let $\ch^n$ be $\mathbb R^{2n}$ equipped with the \emph{Bergman metric} $\bar g$: 
\begin{equation}\label{metrica ch}
\bar g=d\rho^2+\sinh^2(\rho)e_{\cosh^2(\rho)},
\end{equation}
where $\rho$ represents the distance from the origin and $e_{\cosh^2(\rho)}$ is the Berger metric of parameter $\coq$ on $\sss^{2n-1}$. Note that, since the Berger metric changes with the radius, $\bar{g}$ is not given by a warped product.  The Riemann curvature tensor of this metric has the following explicit expression
\begin{equation}\label{curv}
\begin{array}{rcl}
\bar R(X,Y,Z,W)&=&-\bar g(X,Z)\bar g(Y,W)+\bar g(X,W)\bar g(Y,Z)\\&&-\bar g(X,JZ)\bar g(Y,JW)+\bar g(X,JW)\bar g(Y,JZ)\\
&&-2\bar g(X,JY)\bar g(Z,JW),
\end{array}
\end{equation}
where $J$ is the complex structure of $\ch^n$. In our model it coincides with the usual complex structure of $\mathbb R^{2n}$. It follows that the sectional curvature of a plane spanned by two orthonormal vectors $X$ and $Y$ is given by
\begin{equation}\label{sez}
\bar K(X\wedge Y) = -1-3\bar g(X,JY)^2.
\end{equation}
Then $-4\leq\bar K\leq -1$ and it is equal to $-1$ (respectively to $-4$) if and only if $X$ and $JY$ are orthogonal (respectively parallel). This property makes the complex hyperbolic space a \emph{complex space form}. Moreover it is Einstein with $\bar Ric=-2(n+1)\bar g$ and symmetric of rank one, then $\bar \nabla\bar R=0$. 

Here and in the following we are using the convention to put a bar over the symbol for geometric quantity of the fixed ambient manifold $\ch^n$, for example $\bar{\nabla}$ is the Levi-Civita connection of $\bar{g}$.

\subsection{Inverse mean curvature flow}
In \eqref{imcf_def} we defined the inverse mean curvature flow. Since we are considering only closed and mean convex initial data, it is well know that the flow \eqref{imcf_def} has a unique smooth solution, at least for small times. The main geometric quantities for an hypersurface are: the induced metric $g_{ij}$, its inverse $g^{ij}$, the second fundamental form $h_{ij}$, the mean curvature $H=h_{ij}g^{ji}$ and the volume $\left|\mm_t\right|$. They evolve in the following way along the inverse mean curvature flow.

\begin{Lemma}\label{eq_evoluz} Since the ambient space is symmetric the following evolution equations hold:
\begin{itemize}
\item[(1)] $\displaystyle{\frac{\partial g_{ij}}{\partial t} = \frac 2H h_{ij}}$,
\item[(2)] $\displaystyle{\frac{\partial g^{ij}}{\partial t} = -\frac 2H h^{ij}}$,
\item[(3)] $\displaystyle{\frac{\partial H}{\partial t} = \frac{\Delta H}{H^2}-2\frac{\left|\nabla H\right|^2}{H^3}-\frac{\left|A\right|^2}{H}-\frac{\bar{R}ic(\nu,\nu)}{H}}$,
%\item[(4)] $\displaystyle{\frac{\partial h_{ij}}{\partial t} =\frac{\nabla_i\nabla_jH}{H^2}+\frac{h_i^lh_{lj}}{H}-\frac{2}{H^3}\nabla_iH\nabla_jH-\frac{1}{H}\bar R_{i0j0}}$,
\item[(4)] $\displaystyle{\frac{\partial h_{ij}}{\partial t} =\frac{\Delta h_{ij}}{H^2}-\frac{2}{H^3}\nabla_iH\nabla_jH+\frac{\left|A\right|^2}{H^2}h_{ij}-\frac{2}{H}\bar R_{i0j0}+\bar{R}ic(\nu,\nu)\frac{h_{ij}}{H^2}}$
\item[\phantom{(4')}]$\displaystyle{\phantom{\frac{\partial h_{ij}}{\partial t} =}+\frac{1}{H^2}g^{lr}g^{ms}\left(2\bar R_{risj}h_{lm}-\bar R_{rmis}h_{jl}-\bar R_{rmjs}h_{il}\right)}$,
\item[(5)] $\displaystyle{\frac{\partial h_i^j}{\partial t} =\frac{\Delta h_i^j}{H^2}-\frac{2}{H^3}\nabla_iH\nabla_kHg^{kj}+\frac{\left|A\right|^2}{H^2}h_i^j-\frac{2}{H}\bar R_{i0k0}g^{kj}-2\frac{h_i^kh_k^j}{H}}$
\item[\phantom{(5)}]$\displaystyle{\phantom{\frac{\partial h_i^j}{\partial t} =}+\bar{R}ic(\nu,\nu)\frac{h_i^j}{H^2}+\frac{1}{H^2}g^{lr}g^{ms}g^{kj}\left(2\bar R_{risk}h_{lm}-\bar R_{rmis}h_{kl}-\bar R_{rmks}h_{il}\right)}$,
\item[(6)] $\displaystyle{\frac{d\left|\mm_t\right|}{d t}= \left|\mm_t\right|}$,
\item[(7)] $\displaystyle{\frac{\partial\nu}{\partial t}= \frac{\nabla H}{H^2}}$.
\end{itemize}
\end{Lemma}

Here and in the following we are using Einstein convention on repeted indices. Moreover the operation of raising/lowering the indices is done with respect to the induced metric: for example $h_i^j=h_{ik}g^{kj}$. The proof of this Lemma is similar to the computation of the analogous equations for the mean curvature flow which can be found in \cite{Hu}. Note that, integrating equation $(6)$, we have that the inverse mean curvature flow is en expanding flow, precisely $\left|\mm_t\right|=\left|\mm_0\right|e^t$.

\section{Geometry of star-shaped hypersurfaces}\setcounter{equation}{0}\label{sez 3}
All the hypersufaces considered in this paper are star-shaped and $\sss^1$-invariant. In this section we compute the main geometric quantities for a generic star-shaped hypersurface in $\ch^n$ and then we will specify them in case of symmetries. 
Let $F:\sss^{2n-1}\to\ch^n$ be a smooth star-shaped immersion. Up to an isometry of the ambient space, we can consider that it is star-shaped with respect to the origin. Then $F$ is defined by its radial function: there exist a smooth function $\rho:\sss^{2n-1}\to\mathbb R^+$ such that in polar coordinates $\mm=F(\sss^{2n-1})=\left\{(z,\rho(z))\in\ch^n\left|z\in\sss^{2n-1}\right.\right\}$. The hypersurface is said $\sss^1$-invariant if the associated radial function $\rho$ is $\sss^1$-invariant. With the same proof of Lemma 3.1 of \cite{Pi} we can prove the following result. 

\begin{Lemma}
The evolution of an $\sss^1$-invariant hypersurface stays $\sss^1$-invariant during the whole duration of the flow.
\end{Lemma}

Fix any $(Y_1,\dots,Y_{2n-1})$ tangent basis of the sphere $\sss^{2n-1}$, for every $i$ we define $\rho_i=Y_i(\rho)$ and $V_i=F_*Y_i\equiv Y_i+\rho_i\frac{\partial}{\partial \rho}$. Then $(V_1,\dots,V_{2n-1})$ is a tangent basis of $\mm$. The induced metric on $\mm$ is $g=F^*\bar g$, in local coordinates we have
\begin{equation}\label{metrica}
g_{ij}=\rho_i\rho_j+\siq e_{ij},
\end{equation}
where for short $e_{ij}=(e_{\coq})_{ij}$. The inverse of the metric therefore is
\begin{equation}\label{metrica_inv}
g^{ij}=\frac{1}{\siq}\left(e^{ij}-\frac{\rho^i\rho^j}{\siq+|\nabla_e\rho|^2_e}\right),
\end{equation}
where $e^{ij}$ is the inverse of $e_{ij}$, $\rho^i=\rho_ke^{ki}$ and the gradient and the norm of $\rho$ are defined with respect to the metric $e_{\coq}$. In order to simplify the expressions we can fix a function $\varphi=\varphi(\rho)$ such that $\frac{d\varphi}{d\rho}=\frac{1}{\sinh(\rho)}$ and introduce $v=\sqrt{1+\frac{|\nabla_e\rho|^2_{e}}{\siq}}$. Since $\varphi_i=Y_i(\varphi)=\frac{\rho_i}{\sinh(\rho)}$, we get  
\begin{eqnarray*}
g_{ij}&=&\siq(\varphi_i\varphi_j+e_{ij}),\\ g^{ij}&=&\frac{1}{\siq}\left(e^{ij}-\frac{\varphi^i\varphi^j}{v^2}\right),\\v&=&\sqrt{1+|\nabla_e\varphi|^2_e}.
\end{eqnarray*}
A unit normal vector is $$\nu=\frac{1}{v}\left(\parho-\frac{\nabla\rho}{\siq}\right)=\frac{1}{v}\left(\parho-\frac{\nabla\varphi}{\sih}\right).$$
Since the metric of the ambient space is not the same in any direction, it is convenient to use the specific basis tangent to $\sss^{2n-1}$ defined in Lemma \ref{diff_LCC}. In this way we have%to introduce a specific basis tangent to $\sss^{2n-1}$: let us define $(Y_1,\dots,Y_{2n-1})$, where $Y_1=\xi$ is the Hopf vector field and $(Y_2,\dots,Y_{2n-1})$ is $e_{\coq}$-orthonormal and normal to $Y_1$. In this way we have

\begin{equation}\label{berger}
e_{\coq}=\left(\begin{array}{cc}
\coq & 0 \\
0 & id_{2n-2}
\end{array}
\right),
\end{equation}
where $id_{2n-2}$ is the identity matrix of order $2n-2$. The contact form is $\theta(\cdot)=\bar{g}(J\nu,\cdot)$. For an $\sss^1$-invariant star-shaped hypersurface in coordinates we have:
$$
J\nu=\frac{1}{v\sih}\left(\frac{Y_1}{\coh}-\sum_{k\neq 1}\varphi_kJY_k\right),
$$
then
\begin{eqnarray}
\nonumber\theta_i & = & \frac{1}{v\sih}\left(\frac{g_{1i}}{\coh}-\sum_{k\neq 1}\varphi_k\bar{g}(JY_k,V_i)\right)\\
& = & \frac{\sih}{v}\left(\coh\hat{\theta}_i-\sum_{k\neq 1}\varphi_k\delta_{ik^*}\right),\label{theta}
\end{eqnarray}
where
$$
\delta_{ik^*}=\left\{\begin{array}{ccc}
\delta_{i,k+1} & \text{if} & k=2r,\\
-\delta_{i,k-1} & \text{if} & k=2r+1.
\end{array}\right.
$$

Now we want to compute the second fundamental form of $\mm$.  For each $i$ and $j$ let $h_{ij}=-\bar g\left(\bar\nabla_{V_i}V_j,\nu\right)$. Moreover we introduce the following notation: Latin indices $i,j,...$ range from $1$ to $2n-1$ and are related to components tangent to the sphere, the index $0$ represents the radial direction $\parho$ and Greek indices $\alpha, \beta,...$ range from $0$ to $2n-1$. An explicit computation, together to the fact that for every $i$ $\parho(\rho_i)=0$ and $\bar g(\bar\nabla_{Y_i}\parho,\parho)=0$ we get:

$$
h_{ij}=\frac 1v\left(\bar{\Gamma}_{ij}^k\rho_k+\rho_i\rho_k\bar{\Gamma}_{0j}^k+\rho_j\rho_k\bar{\Gamma}_{0i}^k-\bar{\Gamma}_{ij}^0-Y_i(\rho_j)\right).
$$
We have that $\bar{\Gamma}_{ij}^k={\hat{\Gamma}}_{ij}^k$, the Christoffel symbols of the metric $e_{\coq}$, then $$Y_i(\rho_j)-\bar{\Gamma}_{ij}^k\rho_k=\hat{\rho}_{ij},$$
where the ``hat" is in the sense of Notation \ref{notazioni hess}.  For short, let $Y_0$ be $\parho$, then

\begin{eqnarray*}
\bar{\Gamma}_{ij}^0&=&\frac 12\left(Y_i(\bar g_{i\alpha})-Y_{\alpha}(\bar g_{ij})+Y_j(\bar g_{i\alpha})\right)\bar g^{\alpha 0}\\
&= & -\frac 12\parho\left(\bar g_{ij}\right)=-\gamma_i\delta_{ij}=\left\{\begin{array}{ll}
-\sih\coh(\siq+\coq)\delta_{ij} & \text{if } i=1,\\
-\sih\coh\delta_{ij} & \text{if } i\neq 1.  
\end{array}\right.
\end{eqnarray*}

\noindent Finally
\begin{eqnarray*}
\bar{\Gamma}_{i0}^k & = & \frac 12 \left(Y_i\left(\bar g_{0\alpha}\right)-Y_{\alpha}\left(\bar g_{i0}\right)+\parho\left(\bar g_{i\alpha}\right)\right)\bar g^{\alpha k}\\
& =& \frac 12 \parho\left(\bar g_{i\alpha}\right)\bar g^{\alpha k}=\eta_i\delta_{ik}
=\left\{	\begin{array}{ll}
\displaystyle{\frac{\siq+\coq}{\sinh(\rho)\cosh(\rho)}\delta_{ik} }& \text{if } i=1,\\
\displaystyle{\frac{\cosh(\rho)}{\sinh(\rho)}\delta_{ik} }& \text{if } i\neq 1.  
\end{array}\right.
\end{eqnarray*}

\noindent Note that $$\hat{\varphi}_{ij}=\frac{1}{\sinh(\rho)}\hat{\rho}_{ij}-\frac{\cosh(\rho)}{\siq}\rho_i\rho_j\quad\Leftrightarrow\quad\hat{\rho}_{ij}=\sinh(\rho)\hat{\varphi}_{ij}+\sinh(\rho)\cosh(\rho)\varphi_i\varphi_j.$$
Analogous formulas hold for $\rho_{ij}$ and $\varphi_{ij}$ too. Summing up these quantities we get

\begin{eqnarray}
\nonumber h_{ij} & = & \frac{1}{v}\left(-\hat{\rho}_{ij}+\left(\eta_i+\eta_j\right)\rho_i\rho_j+\gamma_i\delta_{ij}\right)\\
 & = & \frac{\sinh(\rho)}{v}\left(-\hat{\varphi}_{ij}+\left(\sinh(\rho)\eta_i+\sih\eta_j-\coh\right)\varphi_i\varphi_j+\frac{\gamma_i}{\sih}\delta_{ij}\right)	\label{sff}
\end{eqnarray}

\noindent Raising the second index we have 
\begin{eqnarray}\label{sff_up}
h_i^k=-\frac{\hat{\varphi}_{ij}\tilde e^{jk}}{v\sih}+\frac{\coh}{v\sih}\delta_i^k+\left\{\begin{array}{ll}
\displaystyle{\frac{\sih}{v\coh}\delta_1^k+\frac{\sih}{v\coh}(\varphi_1)^2\tilde e^{1k}}& \text{ if } i=1;\\
\displaystyle{\frac{\sih}{v\coh}\tilde e_{i1} \tilde e^{1k}}& \text{ if } i\neq 1,
\end{array}\right.
\end{eqnarray}
where $\tilde e^{ij}=\siq g^{ij}=e^{ij}-\frac{\varphi^i\varphi^j}{v^2}$. Taking the trace of \eqref{sff_up} we obtain the mean curvature of $\mm$:

\begin{eqnarray}
\nonumber H & = & h_i^i =-\frac{\hat{\varphi}_{ij}\tilde e^{ji}}{v\sih}+\frac{\sih}{v\coh}+(2n-1)\frac{\coh}{v\sih}\\
\nonumber &&+\frac{\sih}{v\coh}(\varphi_1)^2\tilde{e}^{11}+\frac{\sih}{v\coh}\sum_{i\neq 1}\tilde{e}_{i1}\tilde{e}^{i1}\\
\nonumber & =& -\frac{\hat{\varphi}_{ij}\tilde e^{ji}}{v\sih}+\frac{\sih}{v\coh}+(2n-1)\frac{\coh}{v\sih}\\
\nonumber &&+\frac{\sih}{v\coh}\tilde{e}_{i1}\tilde{e}^{i1}-\sih\coh\frac{\tilde{e}^{11}}{v}\\
\nonumber & = &  -\frac{\hat{\varphi}_{ij}\tilde e^{ji}}{v\sih}+2\frac{\sih}{v\coh}+(2n-1)\frac{\coh}{v\sih}\\
\nonumber &&-\frac{\sih}{v\coh}\left(1-\frac{(\varphi_1)^2}{v^2}\right)\\
& = & -\frac{\hat{\varphi}_{ij}\tilde e^{ji}}{v\sih}+\frac{\hat{H}}{v}+\frac{\sih}{v^3\coh}(\varphi_1)^2,\label{mc}
\end{eqnarray}
where \begin{equation}\label{hat_H}
\hat{H}=\hat{H}(\rho)=(2n-1)\frac{\coh}{\sih}+\frac{\sih}{\coh}=\frac{2n\coq-1}{\sih\coh}.
\end{equation}

If the hypersurface is $\sss^1$-invariant, these expressions can be simplified because, in this case, $\rho_1=\varphi_1=0$. Hence we have:
\begin{eqnarray}
%|\nabla\varphi|^2_{e}&=&|\nabla\varphi|^2_{\sigma};\label{n_phi_s1}\\
h_i^k&=&-\frac{\hat{\varphi}_{ij}\tilde{e}^{jk}}{v\sih}+\frac{\coh}{v\sih}\delta_i^k+\frac{\sih}{v\coh}\delta_i^1\delta_1^k\label{2ff_s1} \\%\frac{1}{v}\left(\begin{array}{cc}
%\frac{\coh}{\sih}+\frac{\sih}{\coh} & 0 \\
%0 & \frac{\sih}{\coh}id_{2n-2}
%\end{array}\right)\\
H & = & -\frac{\hat{\varphi}_{ij}\tilde e^{ji}}{v\sih}+\frac{\hat{H}}{v}.\label{mc_s1}
\end{eqnarray}
The last equation can be also written in a second useful way : let $\tilde{\sigma}^{ij}=\sigma^{ij}-\frac{\varphi^i\varphi^j}{v^2}$, then, by Lemma \ref{hessiani} we have that
$$
\hat{\varphi}_{ij}\tilde{e}^{ij}={\varphi}_{ij}\tilde{\sigma}^{ij},
$$
hence
\begin{equation}\label{mc_s1_sigma}
H=-\frac{{\varphi}_{ij}\tilde {\sigma}^{ji}}{v\sih}+\frac{\hat{H}}{v}.
\end{equation}

\section{The case of geodesic spheres}\label{sez sf geo}

In this section we specify what we found in the previous one in the case of the geodesic spheres and compute their evolution under inverse mean curvature flow. 
Let $\mm_0$ be a geodesic sphere, i.e. a star-shaped hypersurface with constant radial function $\rho=\rho_0$ (of course this function is $\sss^1$-invariant). From \eqref{2ff_s1} we can see that $\mm_0$ has two distinct principal curvatures: $\lambda(\rho)=\coth(\rho)$ with multiplicity $2n-2$ and $\mu(\rho)=\tanh(\rho)+\coth(\rho)=2\coth(2\rho)$ with multiplicity $1$ and eigenvector $\xi$. It follows that 
$H=\hat{H}$. In particular the mean curvature depends only on the radius and then the evolution of $\mm_0$ by inverse mean curvature flow reduces to an ODE: the evolution on $\mm_0$ is a family of geodesic spheres $\mm_t$ of radius $\rho(t)$ satisfying
$$\left\{\begin{array}{l}
\dot{\rho}=\frac 1H=\frac{\sih\coh}{2n\coq-1},\\
\rho(0) = \rho_0.
\end{array}\right.
$$
Trying to solve this ODE we arrive at an implicit value for $\rho(t)$: $$\cosh(\rho(t))\sinh^{2n-1}(\rho(t))=\cosh(\rho_0)\sinh^{2n-1}(\rho_0)e^t.$$ Therefore the solution is defined for any positive time and $\rho(t)=\frac{t}{2n}+o(1)$ as $t\rightarrow\infty$. Moreover $\left|\mm_t\right|=\left|\mm_0\right|e^t$
. Then we get that the rescaled induce metric and the rescaled contact form are:
\begin{eqnarray*}
\tilde g_t&=&\left|\mm_t\right|^{-\frac{1}{n}}g=\frac{\siqt}{\left|\mm_0\right|^{\frac 1n} e^{\frac tn}}e_{\coqt},\\
\tilde{\theta}&=&\left|\mm_t\right|^{-\frac{1}{n}}\theta=\frac{\sinh(\rho(t))\cosh(\rho(t))}{\left|\mm_0\right|e^{\frac 1n}}\hat{\theta}.
\end{eqnarray*}
Obviously $\frac{\siqt}{e^{\frac tn}}\rightarrow \frac 14$,  $\frac{\sinh(\rho(t))\cosh(\rho(t))}{e^{\frac 1n}}\rightarrow \frac 14$ and $\coqt\rightarrow\infty$ as $t\rightarrow\infty$. Hence the contact form converges to a constant multiple of the standard contact form on $\sss^{2n-1}$, then, in particular, the Webster curvature of the limit is constant. Moreover we can see the main new phenomenon of this paper: the rescaled metric does not converge to a Riemannian metric, but to a sub-Riemannian metric defined only on $\mathcal H$. More precisely $\tilde g_t$ converges to a constant multiple of $\sigma_{sR}$. 

The following result is useful to bound the evolution of the radial function in the general case.

\begin{Lemma}\label{confronto}
Consider two concentric geodesic spheres in $\ch^n$ of radius $\rho_1(0)$ and $\rho_2(0)$ respectively. For every $i=1,2$, let $\rho_i(t)$ be the evolution by inverse mean curvature flow of initial datum $\rho_i(0)$, then there is a positive constant $c$ depending only on $n$, $\rho_1(0)$ and $\rho_2(0)$ such that for every time we have
$$
\left|\rho_2(t)-\rho_1(t)\right|\leq c\left|\rho_2(0)-\rho_1(0)\right|.
$$
\end{Lemma}

\proof We can suppose that $\rho_2(0)>\rho_1(0)$ and then this inequality is preserved by the flow. Let us define $\delta=\delta(t)=\rho_2(t)-\rho_1(t)$. The function $\delta$ satisfies
\begin{eqnarray*}
\frac{d\delta}{dt} & = & \frac{1}{(2n-1)\coth(\rho_2)+\tanh(\rho_2)}-\frac{1}{(2n-1)\coth(\rho_1)+\tanh(\rho_1)}
\end{eqnarray*}
Moreover it is easy to see that $\frac{d\rho_1}{dt}>c_1$ with $c_1$ positive constant which depends only on $\rho_1(0)$. Then $\rho_1(t)>c_1 t+\rho_1(0)$. Furthermore trivially $\tanh(\rho_2)>\tanh(\rho_1)>0$ and $\coth(\rho_1)>\coth(\rho_2)>1$. It follows that

\begin{eqnarray*}
\frac{d\delta}{dt}& \leq & \frac{1}{2n-1}\left(\coth(\rho_1)-\coth(\rho_2)\right)\\
& = & \frac{1}{(2n-1)\sinh^2(\tau)}\delta, \qquad \text{for some } \tau\in\left[\rho_1,\rho_2\right]\\
& \leq & \frac{1}{(2n-1)\sinh^2(\rho_1)}\delta\\
& \leq & \frac{1}{(2n-1)\sinh^2(c_1 t+\rho_1(0))}\delta.
\end{eqnarray*}
Integrating we have
\begin{eqnarray*}
\log(\delta(t))-\log(\delta(0))&\leq&\frac{1}{2n-1}\int_0^t\sinh^{-2}(c_1s+\rho_1(0))ds<\infty.\\
%&\leq&\frac{1}{2n-1}\int_0^{\infty}\sinh^{-2}(c_1s+\rho_1(0))ds=c_2<\infty
\end{eqnarray*}
Then the thesis follows.\cvd

From these properties of the geodesic spheres we can deduce some estimates on the evolution of any star-shaped hypersurface. Let $\mm_0$ be defined by the radial function $\rho(0)$, $\rho(t)$ its evolution by inverse mean curvature flow and $\rho_1=\min_{\sss^{2n-1}}\rho(0)$ and $\rho_2=\max_{\sss^{2n-1}}\rho(0)$. Then, with the same notation of the previous Lemma, we have $\rho_1(t)\leq\rho(t)\leq\rho_2(t)$ for any time $t$ the flow is defined. Applying Lemma \ref{confronto} we have that the oscillation of $\rho(t)$ is bounded by a constant which depends only on the initial datum. Below we will show that the flow is defined for any positive time also for any star-shaped $\sss^1$-invariant initial datum. It follows that in any case considered we have $\rho(t)=\frac{t}{2n}+o(1)$ as $t\rightarrow\infty$.

\section{First order estimates}\setcounter{equation}{0}\label{sez 5}
The main result of this section is the proof of part $(1)$ of Theorem \ref{main}. Moreover we will prove also that the mean curvature converges exponentially to that of an horosphere. The main technical result is the following:

\begin{Proposition}\label{gradphi}
There exist a positive constants $c$ such that
$$
|\nabla_{\sigma}\varphi|^2_{\sigma}\leq ce^{-\frac tn}.
$$
%where the gradient and the norm are taken with respect to $\sigma$.
\end{Proposition}

As an immediate geometric consequence we have:

\begin{Corollary}\label{stellato}
The evolution of any star-shaped $\sss^1$-invariant hypersurface stays star-shaped for any time the flow is defined.
\end{Corollary}
\proof An hypersurface is star-shaped if and only if $\parho$ and $\nu$ are never orthogonal in $\ch^n$. This means that there exists a positive constant $c$ such that
$$
\bar g\left(\parho,\nu\right)=\frac 1v\geq \frac 1c \quad\Leftrightarrow\quad v^2\leq c^2
$$
Recalling that $v^2=1+|\nabla_e\varphi|^2_e$, the thesis follows from Proposition \ref{gradphi} noting that $$|\nabla_e\varphi|^2_e=|\nabla_{\sigma}\varphi|^2_{\sigma}$$ holds in case of $\sss^1$ invariance.
\cvd

The proof of Proposition \ref{gradphi} is divided into three steps: first we can prove that $|\nabla_{\sigma}\varphi|^2_{\sigma}$ is bounded, then that  it has an exponential decay and, finally, we show that the right exponent is $\frac 1n$.

\begin{Lemma}\label{gradphi_w}
The following estimate holds:
$$|\nabla_{\sigma}\varphi|^2_{\sigma}\leq\sup_{z\in\sss^{2n-1}}\left|\nabla_{\sigma}\varphi(z,0)\right|^2_{\sigma}.$$
\end{Lemma}
\proof Let us define $\omega=\frac 12 |\nabla_{\sigma}\varphi|^2_{\sigma}=\frac 12\varphi_k\varphi^k$. Note that, as a consequence of the $\sss^1$-invariance, $\varphi_ie^{ik}=\varphi_i\sigma^{ik}$, so we do not distinguish between $\varphi^k$ and $\hat{\varphi}^k$.

We want to compute the evolution equation for $\omega$. We start with the evolution of the radial function:
$$
\frac {1}{Hv}=\frac{d\rho}{dt}=\frac{\partial\rho}{\partial t}+\frac{\partial\rho}{\partial x_i}\frac{\partial x_i}{\partial t}=\frac{\partial\rho}{\partial t}-\frac{\rho^i\rho_i}{Hv\siq}.
$$
Then
\begin{equation}\label{ev_rho}
\frac{\partial\rho}{\partial t}=\frac{1}{Hv}\left(1+\frac{|\nabla_{\sigma}\rho|^2_{\sigma}}{\siq}\right)=\frac vH,
\end{equation}
and so
\begin{equation}\label{ev_phi}
\frac{\partial\varphi}{\partial t}=\frac {1}{\sih}\frac{\partial\rho}{\partial t}=\frac{v}{H\sih}=:\frac 1 F
\end{equation}
holds. From the explicit computation of the mean curvature \eqref{mc_s1_sigma} we have
\begin{eqnarray}\label{F_definizione}
F&=&F(\varphi,\varphi_i,{\varphi}_{ij})=-\frac{{\varphi}_{ij}\tilde {\sigma}^{ij}}{v^2}+\frac{\sih\hat{H}}{v^2}.
%& = & \frac{-\varphi^{(\sigma)}_{ij}\tilde e^{ij}}{v^2}+\frac{\left(\Gamma(e)_{ij}^k-\Gamma(\sigma)_{ij}^k\right)\varphi_k\tilde{e}^{ij}}{v^2}+\frac{\sih\hat{H}}{v^2}
\end{eqnarray}

\noindent Now we want to compute the evolution equation of $\omega$: let $a^{ij}=-\frac{\partial F}{\partial\varphi_{ij}}=\frac{\tilde {\sigma}^{ij}}{v^2}$, $b^i=\frac{\partial F}{\partial\varphi_{i}}$ and, for simplicity of notation, $\nabla=\ns$, then 
\begin{eqnarray*}
\frac{\partial\omega}{\partial t}&=&\varphi ^k\nabla_k\frac{\partial\varphi}{\partial t}\\
& =& -\frac{1}{F^2}\left( -a^{ij}\varphi_{ijk}\varphi^k+b^i\varphi_{ik}\varphi^k+\frac{\partial F}{\partial\varphi}\varphi_k \varphi^k\right)\\
& = &-\frac{1}{F^2}\left( -a^{ij}\varphi_{ijk}\varphi^k +b^i\omega_i+2\frac{\partial F}{\partial\varphi}\omega\right)
\end{eqnarray*}
We can apply the rule for interchanging derivatives:
$$
\varphi_{ijk}=\varphi_{kji}+R^m_{\phantom{m}ijk}\varphi_m,
$$
where this time $R$ is the Riemann curvature tensor of $\sigma$, i.e. $R_{sijk}=\sigma_{sj}\sigma_{ik}-\sigma_{sk}\sigma_{ij}$. Since $a^{ij}$ is symmetric we get:
\begin{eqnarray*}
 -a^{ij}\varphi_{ijk}\varphi^k &= & -a^{ij}\varphi_{kji} \varphi^k-a^{ij}\left(\delta_j^m\sigma_{ik}-\delta_k^m\sigma_{ij}\right)\varphi_m \varphi^k\\
& = & -a^{ij}\omega_{ij}+a^{ij}\varphi_{ik} \varphi_j^k-a^{ij}\varphi_i \varphi_j+2a^i_i\omega.
\end{eqnarray*}
The following equality holds: $$-a^{ij}\varphi_i \varphi_j+2a^i_i\omega=\frac{4(n-1)}{v^2}\omega\geq 0.$$
Moreover, due to the $\sss^1$-invariance, we have 
\begin{equation*}
\frac{\partial}{\partial\rho}\tilde{\sigma}=0,\quad \frac{\partial}{\partial\rho}v=0,
\end{equation*}
hence
\begin{eqnarray}
\nonumber\frac{\partial F}{\partial\varphi} & = & \frac{\partial F}{\partial\rho}\frac{\partial \rho}{\partial\varphi}=\sih\frac{1}{v^2}\frac{\partial}{\partial\rho}\left(\sih\hat{H}\right)\\
  & = & \frac{\siq}{v^2\coq}\left(2n\coq+1\right)>0.\label{stima01}
\end{eqnarray}
Note that $a^{ij}$ is positive definite. Finally we have that 
\begin{eqnarray*}
a^{ij}\varphi_{ik} \varphi_j^k=a^{ij}\sigma^{kl}\varphi_{ik} \varphi_{jl}\geq 0
\end{eqnarray*}
because, as showed in \cite{Di}, if $A$, $B$ and $C$ are symmetric matrices, with $A$ and $B$ positive definite, then $tr(ACBC)\geq 0$.
The thesis follow by the maximum principle. 
\cvd

Now we can use the previous result to bound the mean curvature. In particular we show that the mean convexity is preserved. 

\begin{Lemma}\label{H_bounded}
There exist two positive constants $c_1$ and $c_2$ depending only on $n$ and the initial datum such that for any time the flow is defined
$$
0<c_1\leq H\leq c_2.
$$
\end{Lemma}
\proof From Lemma \ref{eq_evoluz} and the fact that $\left|A\right|^2\geq\frac{1}{2n-1}H^2$ we can compute
$$
\frac{\partial H}{\partial t}\leq \frac{\Delta H}{H^2}-\frac{H}{2n-1}+\frac{2(n+1)}{H}.
$$
By the maximum principle, it is easy to show that $H$ is bounded from above. To prove that H is bounded from below we define $\psi=\frac{v}{\sih H}e^{\frac {t}{2n}}=\frac{1}{F}e^{\frac {t}{2n}}=\frac{\partial\varphi}{\partial t}e^{\frac {t}{2n}}$ and we prove that this function is bounded from above. Preceding like in the proof of Lemma \ref{gradphi_w}:

\begin{eqnarray*}
\frac{\partial \psi}{\partial t} & = & \pat\left(\frac{\partial\varphi}{\partial t}e^{\frac {t}{2n}}\right)\\
&=&-\frac{1}{F^2}\left(-a^{ij}\frac{\partial{\varphi}_{ij}}{\partial t}+b^i\frac{\partial\varphi_i}{\partial t}+\frac{\partial F}{\partial \varphi}\frac{\partial\varphi}{\partial t}\right)e^{\frac {t}{2n}}+\frac{1}{2n}\psi\\
&=& -\frac{1}{F^2}\left(-a^{ij}{\psi}_{ij}+b^i\psi_i+\frac{\partial F}{\partial\varphi}\psi\right)+\frac {1}{2n}\psi
\end{eqnarray*}

\nid From \eqref{stima01} we have that $$\frac{\partial F}{\partial\varphi}\geq 2n\frac{\siq}{v^2},$$ moreover $\frac{1}{F^2}=\psi^2e^{-\frac tn}$. By Lemma \ref{gradphi_w} $v^2$ is bounded. Since the function $\siq e^{-\frac tn}$ is bounded too, we get that

\begin{equation}\label{eq02}
-\frac{1}{F^2}\frac{\partial F}{\partial \varphi}\psi+\frac{1}{2n}\psi\leq-2n\frac{\siq}{v^2} e^{-\frac{t}{n}}\psi^3+\frac{1}{2n}\psi\leq-c\psi^3+\frac{1}{2n}\psi,
\end{equation}
for some positive constant $c$. By the maximum principle we deduce that $\psi$ is bounded. This imply that there is a constant $c>0$ such that $H\geq cv\frac{e^{\frac{t}{2n}}}{\sih}$. The thesis follows since $v\geq 1$ by definition and $\frac{e^{\frac{t}{2n}}}{\sih}$ is bounded.
\cvd

As a consequence we can improve Lemma \ref{gradphi_w}.

\begin{Lemma}\label{gradphi_exp}
There are two positive constants $c$ and $\gamma$ such that:
$$
|\nabla_{\sigma}\varphi|_{\sigma}^2\leq ce^{-\gamma t}.
$$
\end{Lemma}
\proof Since we proved that $H$ is bounded, we have that $\frac{2}{F^2}\frac{\partial F}{\partial\varphi}\geq\gamma>0$ for some $\gamma$. Proceeding like in the proof of Lemma \ref{gradphi_w} we get the thesis.
\cvd

With the help of this estimate, we can refine the result of Lemma \ref{H_bounded}.

\begin{Lemma}\label{Hexp} There is a positive constant $c$ such that:
$$
|H-2n|\leq ce^{-\gamma t},
$$ 
where $\gamma$ is the same of the previous Lemma.
\end{Lemma}
\proof Let $l_i^j$ be equal to $h_i^j+\delta_i^j-\delta_i^1\delta_1^i$. Explicit computations give:
\begin{eqnarray*}
|l|^2 & = & l_i^jl_j^i=|A|^2+2(H-2n)+6(n-1)-2(h_1^1-2);\\
L & = & b_i^i = H-2n+2(2n-1). 
\end{eqnarray*}
Obviously we have 
$$
|l|^2\geq \frac{L^2}{2n-1}.
$$
Moreover, by \eqref{2ff_s1}, $$h_1^1=\frac{1}{v}\left(\frac{\coh}{\sih}+\frac{\sih}{\coh}\right),$$ hence by Lemma \ref{gradphi_exp}
$$
\left|h_1^1-2\right|\leq ce^{-\gamma t}.
$$
By Lemma \ref{eq_evoluz} we get:
\begin{eqnarray}
\nonumber\displaystyle{\left(\frac{\partial}{\partial t}- \frac{\Delta}{H^2}\right)(H-2n)}&\leq & -\frac{1}{H}\left(|A|^2-2(n+1)\right)\\
\nonumber& = & -\frac{1}{H}\left(|l|^2-2(H-2n)-4(2n-1)+2(h^1_1-2)\right)\\
\nonumber& \leq & -\frac 1H\left(\frac{L^2}{2n-1}-2(H-2n)-4(2n-1)\right)+ce^{-\gamma t}\\
\label{eq01}& = & -\frac{H-2n}{(2n-1)H}\left(H+2n-2\right)+ce^{-\gamma t}.
\end{eqnarray}
Hence $H-2n\leq ce^{-\gamma't}$, where $\gamma'\leq\gamma$. Using this information and starting again from \eqref{eq01} we have:
\begin{eqnarray*}
\displaystyle{\left(\frac{\partial}{\partial t}- \frac{\Delta}{H^2}\right)(H-2n)}&\leq & -\left(\frac 1n+ce^{-\gamma't}\right)(H-2n)+ce^{-\gamma t}.
\end{eqnarray*}
Applying again the maximum principle we have the desired estimate from above. On the other side we can consider the evolution of the function $\psi$ defined in Lemma \ref{H_bounded}. Using what we proved so far, the reaction term \eqref{eq02} becomes:
$$
-\frac{1}{F^2}\frac{\partial F}{\partial\varphi}\psi+\frac{1}{2n}\psi \leq \frac{\psi}{2n}\left(1-4n^2\left(1+ce^{-\gamma t}\right)\psi^2\right).
$$
By the maximum principle we have 
$$
\psi\leq\frac {1}{2n}+ce^{-\gamma t},
$$
where the constant $c$ can be different from above. By the definition of $\psi$ we get the thesis.\cvd

Finally we are able to prove the main result of this section.

\noindent \emph{Proof of Proposition \ref{gradphi}.} By the previous Lemma and \eqref{stima01} we have:
$$
\frac{2}{F^2}\frac{\partial F}{\partial\varphi} \geq \frac{4n}{H^2}\geq\frac{1}{n(1+ce^{-\gamma t})}.
$$
With the same computations of the proof of Proposition \ref{gradphi}, and by the maximum principle we have that $|\nabla_{\sigma}\varphi|^2_{\sigma}\leq y$, where $y$ is the solution of the Cauchy problem:
$$
\left\{\begin{array}{rcl}
y'&=&-\frac{1}{n(1+ce^{-\gamma t})}\\
y(0)& = & y_0\geq sup_{t=0}|\nabla_{\sigma}\varphi|^2_{\sigma}.
\end{array}\right.
$$
Then $|\nabla_{\sigma}\varphi|^2_{\sigma}\leq y_0\left(\frac{1+c}{e^{\gamma t}+c}\right)^{\frac{1}{n\gamma}}\leq ce^{-\frac tn}$, for some $c$.
\cvd

We finish this section noting that the proofs of Lemma \ref{Hexp} can be repeated using $\gamma=\frac 1n$, hence we have:

\begin{Proposition}\label{H_lambda_exp}
There is a positive constant $c$ such that:
%\begin{itemize}
%\item[(1)] 
%$$|\lambda_1-2|\leq ce^{-\frac 1n t}.$$
%\item[(2)] 
$$|H-2n|\leq ce^{-\frac 1n t}.$$
%\end{itemize}
\end{Proposition}

\section{Second order estimates}\setcounter{equation}{0}\label{sez 6}

In this section we collect some important results concerning the second order derivative of the radial function. First of all we prove that the principal curvatures of the evolving hypersurface stay bounded. As a consequence we have the long time existence of the flow. In the second part of this section we prove that the Hessian of the radial function is bounded. In view of   \eqref{2ff_s1}, this implies the convergence with exponential speed of the second fundamental form to that one of the horospheres. In this context a new phenomenon appears: in the cases already known in literature, see for example \cite{Ge1,Ge2,Sc,Zh}, we always have the same speed of convergence for any initial datum. Here, in the general case, we can expect just half of the speed of the convergence of, for example, the evolution of a geodesic sphere. See Remark \ref{remark} below for more details.

%We already know that the principal curvature associate to the special direction tangent to the fibers of the action of $\sss^1$ converges everywhere to $2$. In this section we want to study the behaviour of all other principal curvatures. First of all we  prove part $(2)$ of Theorem \ref{main}, showing that the principal curvatures of the evolving hypersurface are uniformly bounded. %As a consequence we are able to improve the results of the previous sections: we will prove the mean curvature tends exponentially fast to $2n$ (the value of the mean curvature of an horosphere) and that $\gradphi_{\sigma}\leq ce^{-\frac tn}$.
 
\begin{Lemma}\label{eq_M}
Let us define the tensor $M_i^{j}=Hh_i^{j}$, then 
\begin{eqnarray*}
\frac{\partial M_i^j}{\partial t} &=& \frac{\Delta M_i^j}{H^2}-\frac{2}{H^3}\left\langle \nabla H,\nabla M_i^j\right\rangle-\frac{2}{H^3}\nabla_iH\nabla_kHg^{kj}-2\frac{M_i^kM_k^j}{H^2}-2\bar R_{i0k0}g^{kj}\\
 &&+\frac{1}{H^2}g^{qs}g^{kj}\left(2\bar R_{pisk}M_q^p-\bar R_{pqis}M_k^p-\bar R_{pqks}M_i^p\right).
\end{eqnarray*} 
\end{Lemma}

\proof Some simple computations give that  
\begin{eqnarray*}
\nabla M_i^j&=&h_i^j\nabla H+H\nabla h_i^j\\		
\Delta M_i^j &= &h_i^j\Delta H+H\Delta h_i^j+2\left\langle \nabla H,\nabla h_i^j\right\rangle\\
& = & h_i^j\Delta H+H\Delta h_i^j+\frac{2}{H}\left\langle \nabla H,\nabla M_i^j\right\rangle-2\frac{h_i^j}{H}\left|\nabla H\right|^2
\end{eqnarray*}

\nid Using the evolution equation for $H$ and $h_i^j$ in Lemma \ref{eq_evoluz} we get the thesis.
\cvd

\begin{Corollary}\label{2ff_bound}
The principal curvatures of the hypersurface evolving by inverse mean curvature flow are uniformly bounded.
\end{Corollary}
\proof Since by Lemma \ref{H_bounded} $H$ is bounded (in particular from below), it is sufficient to prove that the principal curvatures are bounded from above. Let $m=2n-1$ and $\mu_1\leq\dots\leq\mu_{m}$ the eigenvalues of $M_i^j$. The trace of $M_i^j$ is $\sum_i\mu_i=H^2>0$, then $\mu_{m}>0$ everywhere. To conclude the proof we want to show that $\mu_{m}$ is bounded from above. Fix any time $T^*$ strictly smaller than the maximal time $T$. We can find a point $(x_0,t_0)$ where $\mu_{m}$ reaches its maximum in $\sss^{2n-1}\times\left[0,T^*\right]$. At this point we can fix an orthonormal basis which diagonalizes $M_i^j$, then we can say that at this point $\mu_{m}$ satisfies the same equation for $M_{m}^{m}$ found in Lemma \ref{eq_M}. The following estimates hold:

\begin{eqnarray*}
-\frac{2}{H^2}\nabla_{m}H\nabla_kHg^{km} & = & -\frac{2}{H^2}(\nabla_{m}H)^2\leq 0,\\
-\bar R_{m 0 j 0}g^{jm} & = & -\bar R_{m0m0}\leq 4,\\
2g^{qs}g^{km}\bar R_{pmsk}M_q^p&=&\sum_q2\mu_q\bar R_{qmqm},\\
g^{qs}g^{km}\bar R_{pqms}M_k^p&=&g^{qs}g^{km}\bar R_{pqks}M_m^p\\
&=&\sum_q\mu_m\bar R_{qmqm}.
\end{eqnarray*}

\nid It follows that
\begin{eqnarray*}
g^{qs}g^{km}\left(2\bar R_{pmsk}M_q^p-\bar R_{pqms}M_k^p-\bar R_{pqks}M_m^p\right)& =& 2\sum_q\bar R_{qmqm}\left(\mu_q-\mu_m\right)\\
& \leq & -8\sum_q\left(\mu_q-\mu_m\right)\\
 & = & 8(2n-1)\mu_m-8H^2.
\end{eqnarray*}
Putting together these computations we have that at $(x_0,t_0)$ 

$$
\frac{\partial\mu_m}{\partial t} \leq \frac{\Delta\mu_m}{H^2}-\frac{2}{H^3}\left\langle \nabla H,\nabla\mu_m\right\rangle-\frac{2\mu_m}{H^2}\left(\mu_m-4(2n-1)\right).
$$
If $t_0\neq 0$ we can deduce that in $(x_0,t_0)$
$$
0\leq \frac{\partial\mu_m}{\partial t}-\frac{\Delta\mu_m}{H^2}+\frac{2}{H^3}\left\langle \nabla H,\nabla\mu_m\right\rangle\leq-\frac{2\mu_m}{H^2}\left(\mu_m-4(2n-1)\right).
$$
Since $\mu_m$ is positive, it follows that $\mu_m\leq 4(2n-1)$. Hence $\mu_m$ reaches its maximum at time $t_0=0$ or it is bounded by a constant independent on the choice of $T^*$.
 \cvd

It follows the uniform parabolicity of equation \eqref{ev_rho} and an uniform $C^2$-estimate for the function $\rho(\cdot, t)$. Arguing as in chapter 2.6 of \cite{Ge1}, we can  apply the $C^{2,\alpha}$ estimates of \cite{Kr} to conclude that the solution of the flow is defined for any positive time and it is smooth, since the initial datum is smooth. 

In the sequel we can perform a better analysis giving an uniform estimate on the Hessian of $\rho$. We prefer to work with the auxiliary function $\varphi$.

\begin{Proposition}\label{hess_phi_exp}
There is a positive constant $c$ such that
$$
\left|\nabla_{\sigma}^2\varphi\right|^2_{\sigma}\leq ce^{-\frac tn}.
$$
\end{Proposition}
\proof The proof is similar to that of Lemma \ref{gradphi_w}, but we will use more than once the interchanging rule. Let us define this time $\omega=\frac 12\left|\nabla_{\sigma}^2\varphi\right|^2_{\sigma}$, then, recalling Notations \ref{notazioni hess} we have:%\textbf{riv notazione ordine indici}
\begin{eqnarray*}
\frac{\partial}{\partial t}\omega & = & \varphi^{rs} \nabla_r\nabla_s\frac{\partial}{\partial t}\varphi\\
& = & \frac{1}{F^2}\varphi^{rs}\left(a^{ij}\varphi_{ijrs}-b^i\varphi_{irs}-\frac{\partial F}{\partial \varphi}\varphi_{rs}\right)\\
& = &\frac{1}{F^2}\left(a^{ij}\varphi^{rs} \varphi_{ijrs} -b^i\varphi^{rs}\varphi_{irs} -2\frac{\partial F}{\partial \varphi}\omega\right).
\end{eqnarray*}

\nid Applying twice the interchanging rule we have:

$$
\begin{array}{l}
a^{ij}{\varphi}^{rs}{\varphi}_{rsij}\\
\qquad = a^{ij}{\varphi}^{rs}\left({\varphi}_{isrj}+(\delta^m_i\sigma_{sr}-\delta^m_r\sigma_{si}){\varphi}_{mj}\right)\\
\qquad = a^{ij}{\varphi}_{rs}\left[\nabla_i\left({\varphi}_{jrs}+(\delta^m_j\sigma_{sr}-\delta^m_s\sigma_{rj}){\varphi}_m\right)\right]\\
\qquad\phantom{=} +a^{ij}{\varphi}_{ij}\Delta\varphi-a^{ij}{\varphi}^r_i {\varphi}_{rj}\\
\qquad = a^{ij}{\varphi}_{rs} {\varphi}_{ijrs}+2a^{ij}{\varphi}_{ij}\Delta\varphi-2a^{ij}{\varphi}^r_i {\varphi}_{rj}\\
\qquad = a^{ij}{\omega}_{ij}+2a^{ij}{\varphi}_{ij}\Delta\varphi-2a^{ij}{\varphi}^r_i {\varphi}_{rj}\\
\qquad\phantom{=}-a^{ij}{\varphi}_i^{rs} {\varphi}_{jrs};
\end{array}
$$

%$$
%\begin{array}{l}
%a^{ij}D^rD^s\varphi D_rD_sD_iD_j\varphi\\
%\qquad = a^{ij}D^rD^s\varphi\left(D_iD_sD_rD_j\varphi+(\delta^m_i\sigma_{sr}-\delta^m_r\sigma_{si})D_mD_j\varphi\right)\\
%\qquad = a^{ij}D^rD^s\varphi\left[D_i\left(D_jD_rD_s\varphi+(\delta^m_j\sigma_{sr}-\delta^m_s\sigma_{rj})D_m\varphi\right)\right]\\
%\qquad\phantom{=} +a^{ij}D_iD_j\varphi\Delta\varphi-a^{ij}D^rD_i\varphi D_rD_j\varphi\\
%\qquad = a^{ij}D^rD^s\varphi D_iD_jD_sD_r\varphi+2a^{ij}D_iD_j\varphi\Delta\varphi-2a^{ij}D^rD_i\varphi D_rD_j\varphi\\
%\qquad = a^{ij}D_iD_j\omega+2a^{ij}D_iD_j\varphi\Delta\varphi-2a^{ij}D^rD_i\varphi D_rD_j\varphi\\
%\qquad\phantom{=}-a^{ij}D_iD^rD^s\varphi D_jD_rD_s\varphi;
%\end{array}
%$$
\nid while
$$
\begin{array}{l}
-b^i\varphi ^{rs} \varphi_{irs}\\
\qquad = -b^i\varphi^{rs}\left(\varphi_{rsi}+(\delta^m_i\sigma_{rs}-\delta^m_r\sigma_{ir})\varphi_m\right)\\
\qquad = -b^i\varphi^{rs} \varphi_{rsi}-b^i\varphi_i\Delta\varphi+b^i\varphi_i^r \varphi_r.
\end{array}
$$
Summing up these quantities we have:
$$
\frac{\partial}{\partial t}\omega = \frac{a^{ij}}{F^2}\omega_{ij}-\frac{b^i}{F^2}\omega_i-\frac{2}{F^2}\frac{\partial F}{\partial\varphi}\omega+\frac{R}{F^2},
$$
where the remainder term is
\begin{eqnarray*}
R & = & 2a^{ij}\varphi_{ij}\Delta\varphi-2a^{ij}\varphi^r_j \varphi_{ri}-a^{ij}\varphi_i^{rs} \varphi_{jrs}\\
 & & -b^i\varphi_i\Delta\varphi+b^i\varphi^r_i \varphi_r.
\end{eqnarray*}
From \eqref{stima01} and Proposition \ref{Hexp} we get
$$
-\frac{2}{F^2}\frac{\partial F}{\partial\varphi} = -\frac{2}{H^2}\left(2n+\frac{1}{\coq}\right)\leq-\frac 1n + ce^{-\frac tn}.
$$
Moreover we claim that 
$$
R\leq c(\omega+1).
$$
In fact we notice that, for $t$ big enough, $a^{ij}\geq \frac 12 \sigma^{ij}$, hence
\begin{eqnarray*}
-a^{ij}\varphi_i^{rs} \varphi_{jrs}&\leq& -\frac 12|\nabla^3_{\sigma}\varphi|^2_{\sigma}\leq 0,\\
-2a^{ij}\varphi_i^r \varphi_{rj}  & \leq& c\omega.
\end{eqnarray*}

\noindent Obviously $(\Delta\varphi)^2\leq 2(2n-1)\omega$ and, by \eqref{mc_s1_sigma}, 
$$
a^{ij}\varphi_{ij}=\frac{\sih}{v}\left(\frac{\hat{H}}{v}-H\right),
$$
therefore
$$
2a^{ij}\varphi_{ij}\Delta\varphi\leq\left(a^{ij}\varphi_{ij}\right)^2+\left(\Delta\varphi\right)^2\leq c(\omega+1).
$$
An explicit computation shows that
$$
b^i=-\frac{2F}{v^2}\varphi^i+\frac{{\varphi}_{rs}}{v^4}\left(\varphi^r\delta^{si}+\varphi^s\delta^{ri}-2\frac{\varphi^r\varphi^s\varphi^i}{v^2}\right).
$$
It follows that
\begin{eqnarray*}
|b^i|^2_{\sigma} & \leq & c(\omega+1),\\
%b^iD_i\varphi & = & \frac{2F}{v^2}|\nabla\varphi|^2_{\sigma}+\frac{2}{v^6}\hat{\varphi}_{rs}\varphi^r\varphi^s.
\end{eqnarray*}
Using these results we can finally estimate
\begin{eqnarray*}
-b^i\varphi_i\Delta\varphi & \leq& \frac 12|b^i|^2_{\sigma}|\nabla_{\sigma}\varphi|^2_{\sigma}+\frac 12(\Delta\varphi)^2\leq c(\omega+1),\\
b^i\varphi^r\varphi_{ri} & \leq &  c(|b^i|^2_{\sigma}+|\nabla\varphi|^2_{\sigma}+\omega))\leq c(\omega+1).
\end{eqnarray*}
Summing up what we found we have that
$$
\frac{\partial}{\partial t}\omega \leq \frac{a^{ij}}{F^2}{\omega}_{ij}-\frac{b^i}{F^2}\omega_i+\left(-\frac 1n+ce^{-\frac tn}\right)\omega+ce^{-\frac tn}.
$$
The thesis follows by the maximum principle.\cvd

\begin{Remark}
In \eqref{ev_rho} we computed the scalar flow satisfied by the radial function. The flow \eqref{imcf_def} is defined at least as the corresponding flow for $\rho$. We computed two different expressions for the mean curvature - \eqref{mc_s1} and \eqref{mc_s1_sigma}- each one define, formally, a different scalar flow. However, as consequence of Lemma \ref{hessiani}, in the special case of $\sss^1$-invariance the two flows coincide, so, from a technical point of view, the choice of study the derivatives with respect to $\sigma$, and not $e_{\coq}$ is coherent.
\end{Remark}

A consequence of this Proposition is the convergence of the second fundamental form to that one of an horosphere.

\begin{Corollary}
There is a positive constant $c$ such that
$$
\left|h_i^k-\delta_i^k-\delta_i^1\delta_1^k\right|^2\leq c e^{-\frac tn}.
$$
Moreover on the horizontal distribution we have a faster convergence: taking the sum over $i,k\neq 1$
$$
\left|h_i^k-\delta_i^k\right|^2\leq c e^{-\frac {2t}{n}}.
$$
\end{Corollary}
\proof
By the expression of the second fundamental form \eqref{2ff_s1}, of the mean curvature \eqref{mc_s1}  and Lemma \ref{hessiani}  we have:
\begin{eqnarray}
\nonumber\left|h_i^k-\delta_i^k-\delta_i^1\delta_1^k\right|^2&= & \frac{\hat{\varphi_{ij}}\hat{\varphi}_{kr}\tilde{e}^{jk}\nonumber\tilde{e}^{ri}}{v^2\siq}-\frac{2\hat{\varphi}_{ij}\tilde{e}^{ji}}{v\sih}\left(\frac{\coh}{v\sih}-1\right)\\
\nonumber&&+(2n-1)\left(\frac{\coh}{v\sih}-1\right)^2+\left(\frac{\sih}{v\coh}-1\right)^2\\
\nonumber&&+2\left(\frac{\coh}{v\sih}-1\right)\left(\frac{\sih}{v\coh}-1\right)\\
\nonumber&=& \frac{{\varphi_{ij}}{\varphi}_{kr}\tilde{\sigma}^{jk}\tilde{\sigma}^{ri}}{v^2\siq}+2\frac{|\nabla_{\sigma}\varphi|^2_{\sigma}}{v^2}+2\left(H-\frac{\hat{H}}{v}\right)\left(\frac{\coh}{v\sih}-1\right)\\
\label{eqX}&&+(2n-1)\left(\frac{\coh}{v\sih}-1\right)^2+\left(\frac{\sih}{v\coh}-1\right)^2\\
\nonumber&&+2\left(\frac{\coh}{v\sih}-1\right)\left(\frac{\sih}{v\coh}-1\right)\\
\end{eqnarray}
By Proposition \ref{gradphi}, Proposition \ref{H_lambda_exp} and Proposition \ref{hess_phi_exp} we get that all the terms appearing in the last equality can be bounded by $ce^{-\frac {2t}{n}}$, except the ``bad" term $2\frac{|\nabla_{\sigma}\varphi|^2_{\sigma}}{v^2}$ that is just smaller than $ce^{-\frac tn}$. Hence the first estimate is proven.
\nid Finally, if we restrict our attention to the horizontal distribution, the same computations give:
\begin{eqnarray*}
\left|h_i^k-\delta_i^k\right|^2 &=& \frac{{\varphi_{ij}}{\varphi}_{kr}\tilde{\sigma}^{jk}\tilde{\sigma}^{ri}}{v^2\siq}+2\left(H-\frac{\hat{H}}{v}\right)\left(\frac{\coh}{v\sih}-1\right)\\
&&+(2n-2)\left(\frac{\coh}{v\sih}-1\right)^2\leq ce^{-\frac{2t}{n}}.
\end{eqnarray*}

\cvd

\begin{Remark}\label{remark}
For some initial data we are able to find the ``fast" convergence of the second fundamental form in the whole tangent space, because $|\nabla_{\sigma}\rho|^2_{\sigma}$  converges to (or is, like in the trivial case of geodesic spheres) zero. In general the estimate found in the previous Corollary is optimal. In fact, for the examples that we will discuss in the last section, we have that the gradient of the radial function is just bounded and cannot converges to zero. Hence Proposition \ref{gradphi} cannot be improved and the ``bad" term in equation \eqref{eqX} decays slower than all others. Moreover Proposition \ref{H_lambda_exp} says that we don't see this difference at the level of the mean curvature and we find the optimal speed in any case. The reason is that, like shown in Lemma \ref{hessiani}, if $\varphi$ is an $\sss^1$-invariant function, then $\Delta_e\varphi=\Delta_{\sigma}\varphi$. 
\end{Remark}

\section{Higher order estimates}\setcounter{equation}{0}\label{sez 7}

Following the same procedure of the previous section, we can show that the spatial derivative of any order of $\varphi$ has an exponential decay. 

\begin{Notation}
\begin{itemize}
\item[1)] In order to avoid confusion with the meaning of the indices, in this section we use the following notation: capital letters count the number of derivations (for example $\nabla^K\varphi$ is the $K$-th derivative of the function $\varphi$), while lowercase letters are indices representing a direction.
\item[2)] In the next proof we will use also a well established notation: given two tensors $S$ and $T$, we write $S\ast T$ for any linear combination formed by contraction on $S$ and $T$ by $\sigma$.
\end{itemize}
\end{Notation}

\begin{Proposition}\label{higher order}
For any integer $K$, there is a positive constant $c$ which depends only on $K$, $n$, and $\mm_0$, such that
$$
|\nabla^K_{\sigma}\varphi|^2_{\sigma}\leq ce^{-\frac tn}.
$$
\end{Proposition}
\proof The proof follows the same strategy of Proposition \ref{hess_phi_exp}. Fix $K$ and let $\omega=\frac 12 |\nabla^K_{\sigma}\varphi|^2_{\sigma}$. Applying a finite number of times the interchanging rule for derivative we have:
\begin{eqnarray*}
\frac{\partial}{\partial t}\omega & = & \nabla^{i_1}\dots\nabla^{i_K}\nabla_{i_1}\dots\nabla_{i_K}\frac{\partial}{\partial t}\varphi\\
& = & \frac{a^{ij}}{F^2}\omega_{ij}-\frac{b^i}{F^2}\omega_i-\frac{2}{F^2}\frac{\partial F}{\partial\rho}\omega\\
&&-\frac{2a^{ij}}{F^2}\nabla_i\nabla_{i_1}\dots\nabla_{i_K}\varphi\nabla_j\nabla^{i_1}\dots\nabla^{i_K}\varphi\\
&&+\frac{1}{F^2}\left(a\ast\nabla^K\varphi\ast\nabla^K\varphi\right)+\frac{1}{F^2}\left(b\ast\nabla^K\varphi\ast\nabla^{K-1}\varphi\right),\\
%&&-\frac{2}{F^2}\frac{\partial}{\partial\rho}F.
\end{eqnarray*}
Arguing as in the proof of Proposition \ref{hess_phi_exp}, and supposing by induction that the thesis holds for $K-1$, we can show that in this case too we have:
$$
\frac{\partial}{\partial t}\omega \leq \frac{a^{ij}}{F^2}{\omega}_{ij}-\frac{b^i}{F^2}\omega_i+\left(-\frac 1n+ce^{-\frac tn}\right)\omega+ce^{-\frac tn}.
$$
The thesis follows once again as an application of the maximum principle.
\cvd

An immediate consequence of Proposition \ref{higher order} is the following.

\begin{Corollary}\label{higher order rho}
For any integer $K$ there is a positive constant $c$ which depends only on $K$, $n$, and $\mm_0$, such that
$$
|\nabla^K_{\sigma}\rho|^2_{\sigma}\leq c.
$$
\end{Corollary}

\section{Convergence of the rescaled metric and contact form}\setcounter{equation}{0}\label{sez 8}
In this section we prove part \emph{(3)}  of Theorem \ref{main} studying the limit of the rescaled metric and of the rescaled contact form. 

\begin{Theorem}
There is a smooth $\sss^1$-invariant function $f$ such that the metric $\tilde g_t=\left|\mm_t\right|^{-\frac 1n}g_t$ converges to a sub-Riemannian metric $\tilde g_{\infty}=e^{2f}\sigma_{sR}$ and the contact form $\tilde{\theta}=\left|\mm_t\right|^{-\frac 1n}\theta_t$ converges to the conformal multiple $\tilde{\theta}_{\infty}=e^{2f}\hat{\theta}$.
\end{Theorem}
\proof For every time $t$, let $\tilde\rho_t$ be the radius of a geodesic sphere $B_{\tilde\rho_t}$ such that $\left|\mm_t\right|=\left|B_{\tilde\rho_t}\right|$ and define $f_t=\rho(x,t)-\tilde\rho_t$. The mean curvature of $B_{\tilde\rho_t}$ is $\tilde{H}=\hat{H}(\tilde{\rho})$, then $\frac{\partial\tilde{\rho}}{\partial t}=\tilde{H}^{-1}$ and $\tilde\rho=\frac{t}{2n}+o(1)$ as $t\rightarrow\infty$. 

\nid We recall that $g_{ij}=\siq\left(\varphi_i\varphi_j+e_{ij}\right)$. Obviously $e_{\coq}\rightarrow\sigma_{sR}$ as $t\rightarrow\infty$. By Proposition \ref{gradphi} we have that each $\varphi_i$ is going to zero, hence
\begin{eqnarray*}
\lim_{t\rightarrow\infty}\tilde g_t &=&\frac{1}{|\mm_0|^{\frac 1n}}\left(\lim_{t\rightarrow\infty}{\siq}{e^{-\frac{t}{n}}}\right)\sigma_{sR}\\
&=&\gamma\left(\lim_{t\rightarrow\infty}e^{2f_t}\right)\sigma_{sR},
\end{eqnarray*}
for some positive constant $\gamma$.

\nid In analogous way, recalling the expression of the contact form \eqref{theta}, we get
\begin{eqnarray*}
\lim_{t\rightarrow\infty}\tilde{\theta}_t &=&\frac{1}{|\mm_0|^{\frac 1n}}\left(\lim_{t\rightarrow\infty}{\sih\coh}{e^{-\frac{t}{n}}}\right)\hat{\theta}\\
%& = & sss\\
&=&\gamma\left(\lim_{t\rightarrow\infty}e^{2f_t}\right)\hat{\theta},
\end{eqnarray*}
where $\gamma$ is the same constant appearing in the limit of $\tilde{g}_t$.

\nid The thesis follows if we can prove that $f_t$ converges, as $t$ goes to infinity, to a smooth functions $f_{\infty}$, then $f$ will be a multiple of such $f_{\infty}$. Since for any integer $K$, $|\nabla_{\sigma}^Kf_t|^2_{\sigma}=|\nabla_{\sigma}^K\rho|^2_{\sigma}$ and it is uniformly bounded by Corollary \ref{higher order rho}, we have to show that there exist a positive constant $c$ such that $\left|\frac{\partial f_t}{\partial t}\right|\leq ce^{-\frac tn}$. By \eqref{ev_rho}, we have that
$$
\frac{\partial f_t}{\partial t}=\frac{\partial \rho}{\partial t}-\frac{d\tilde\rho_t}{d t}=\frac{v}{H}-\frac{1}{\tilde H}
$$
Using triangle inequality we have
$$
\left|\frac{\partial f_t}{\partial t}\right|\leq\frac{1}{H}\left|v-1\right|+\frac{1}{H\tilde H}\left(\left|H-2n\right|+\left|\tilde H-2n\right|\right).
$$
Since we know that $H$, $\tilde H$ and $v$ are bounded and positive, the desired estimate follows by Proposition \ref{gradphi} and Proposition \ref{H_lambda_exp}.
%\textbf{dire convergenza c infinito su H?}
\cvd

\section{The curvature of the limit metric}\label{sez 9}

In this section we conclude the proof of Theorem \ref{main} showing that $\tilde{\theta}_{\infty}=e^{2f}\hat{\theta}$ (or equivalently $\tilde g_{\infty}=e^{2f}\sigma_{sR}$) does not have necessarily constant Webster scalar curvature. We recall that, as shown in Lemma \ref{caratterizzazione_curv}, since we are considering only $\sss^1$-invariant hypersurfaces, Jerison and Lee's formula \eqref{JL} can be simplified and we have to check just that some of such functions $f$ are not constant.

The construction of the required hypersurface is inspired to the solution of the analogous problem in the real hyperbolic space \cite{HW}. However it is well known that in $\mathbb{CH}^n$ there are no totally umbilical hypersurfaces, so the trace-free part of the second fundamental form cannot have the same strong meaning that it has in the case of hyperbolic space: see  Propositions 3 and 5 of \cite{HW} . So we defined a Brown-York type quantity on hypersurfaces which gives a measure of how the hypersurface is far from being a geodesic sphere. For any star-shaped hypersuface $\mm$ we define 
\begin{equation}\label{Q_def}
Q(\mm)=|\mm|^{-1+\frac {1}{n}}\int_{\mm}\left(H-{\hat H}\right)d\mu,
\end{equation}
where $\hat H$ was defined in \eqref{hat_H}. $Q$ is not a true measure, because we do not know its sign: it is trivially zero when $\mm$ is a geodesic sphere, but in general it is not true the opposite. One of the main property of $Q$ is the following.

%Note that we do not know the sign of this function, but, we are sure that $Q(\mm)$ is bounded along the flow 

\begin{Proposition}\label{LQ}
Let $\widetilde{\mm}^{\tau}$ be a family of hypersurfaces in $\ch^n$ that are radial graph of the functions $\tilde{\rho}(z,\tau)=\tau+f(z)+o(1)$, for some fixed $\sss^1$-invariant function $f:\sss^{2n-1}\rightarrow\mathbb R$. Then
$$
\lim_{\tau\rightarrow\infty}Q(\widetilde{\mm}^{\tau})=\left(\int_{\sss^{2n-1}}e^{2nf}d\sigma\right)^{-1+\frac 1n}\int_{\sss^{2n-1}}e^{2nf}\left(e^{-f}\Delta_{\sigma}e^{-f}-n|\nabla_{\sigma} e^{-f}|_{\sigma}^2\right) d\sigma.
$$
Moreover if $\lim_{\tau\rightarrow\infty}Q(\widetilde{\mm}^{\tau})\neq 0$, then $e^{2f}\sigma_{sR}$ - the limit of the rescaled metric on $\widetilde{\mm}^{\tau}$ - does not have constant Webster curvature.
\end{Proposition}

\proof First of all, note that, since we are considering only $\sss^1$-invariant hypersurfaces, the contribution of the special direction is ruled out and then we can consider the usual Rieamannian Lapacian, the gradient and the volume form associated to $\sigma$ even if we know that the limit metric is sub-Riemannian.
From the expression of the mean curvature of a star-shaped hypersurface \eqref{mc_s1_sigma} we have that
$$
H-{\hat H}=-\frac{\varphi_{ij}\tilde {\sigma}^{ij}}{v\sih}+\hat{H}\left(\frac{1}{v}-1\right)=-\frac{\varphi_{ij}\tilde {\sigma}^{ij}}{v\sih}-\frac{\hat{H}}{v(v+1)}|\nabla_{\sigma}\varphi|_{\sigma}^2.
$$

Since $\tilde{\rho}_1=f_1=0$ for every $j$, we can compute: 
\begin{equation}\label{stime tau}
\left\{\begin{array}{rcl}
v^2 & = & 1+\frac{1}{\sinh^2(\tilde{\rho})}|\nabla_{\sigma}\tilde{\rho}|_{\sigma}^2\\
& =& 1+O(e^{-{\tau}});\\
\hat{H} & = & 2n+O(e^{-{\tau}});\\
\tilde{\sigma}^{ij} &=& \sigma^{ij}-\frac{\tilde{\rho}^i\tilde{\rho}^j}{v^2\siq}\\
& =& \sigma^{ij}+O(e^{-{\tau}});\\
\varphi_i & = & \frac{\tilde{\rho}_i}{\sinh(\tilde{\rho})} = \frac{e^f}{\sinh(\tilde{\rho})}\nabla_ie^{-f};\\
\varphi_{ij} & = & \frac{1}{\sinh(\tilde{\rho})}\left(\tilde{\rho}_{ij}-\frac{\cosh(\tilde{\rho})}{\sinh(\tilde{\rho})}\tilde{\rho}_i\tilde{\rho}_j\right)\\
& =& \frac{1}{\sinh(\tilde{\rho})}\left(f_{ij}-f_if_j+o(1)\right)\\
 & =& -\frac{e^f}{\sinh(\tilde{\rho})}\left(\nabla^2_{ij}e^{-f}+o(1)\right).
\end{array}\right.
\end{equation}

\nid It follows that
\begin{eqnarray*}
\lim_{\tau\rightarrow\infty}Q(\widetilde{\mm}^{\tau})&=& \lim_{\tau\rightarrow\infty}\left[\left(\int\sinh^{2n-1}(\tilde{\rho})\cosh(\tilde{\rho})d\mu_{\tau}\right)^{-1+\frac 1n}\right.\\
&&\phantom{aaa}\left. * \int\left(\sinh^{2n-1}(\tilde{\rho})\cosh(\tilde{\rho})\left(-\frac{\varphi_{ij}\tilde{\sigma}^{ij}}{v\sinh(\tilde{\rho})}-\frac{\hat{H}}{v(v+1)}(v^2-1)\right)\right)d\mu_{\tau}\right]\\
& = & \left(\int_{\sss^{2n-1}}e^{2nf}d\sigma\right)^{-1+\frac 1n}\int_{\sss^{2n-1}}e^{2nf}\left(e^{-f}\Delta_{\sigma}e^{-f}-n|\nabla_{\sigma} e^{-f}|_{\sigma}^2\right) d\sigma.
\end{eqnarray*}

\nid This formula shows that if $\lim_{\tau\rightarrow\infty}Q(\widetilde{\mm}^{\tau})\neq 0$, then $e^{-f}\Delta_{\sigma}e^{-f}-n|\nabla_{\sigma} e^{-f}|^2_{\sigma}\neq 0$ and so $f$ cannot be constant. Lemma \ref{caratterizzazione_curv} tells us that the limit metric $e^{2f}\sigma_{sR}$ does not have constant Webster curvature. Finally note that the opposite is not true because $e^{-f}\Delta_{\sigma}e^{-f}-n|\nabla_{\sigma} e^{-f}|^2_{\sigma}$ does not have necessarily a sign.
\cvd

If we compare $Q$ with the modified Hawking mass studied for the real hyperbolic case in \cite{HW}, $Q$ has the disadvantage that it works only with $\sss^1$-invariant data and it does not characterize the constant curvature limit. However Proposition \ref{LQ} suggests that the study of the asymptotic behaviour of $Q$ is enough to find a family of initial data such that the limit of the rescaled metric does not have constant Webster curvature. In order to complete this goal we need to study the evolution equation of $Q$.

\begin{Lemma}\label{eq_evoluz_2}
%\begin{itemize}
%\item[(1)] 
%\item[(2)] If $\mm_0$ is $\sss^1$-invariant we have $$\displaystyle{\frac{\partial }{\partial t}\frac{1}{v}=\frac{1}{H^2}\left\langle\nabla H,\partial \rho^T\right\rangle+\frac{\coh}{H\sih}\frac{v^2-1}{v^2}},$$
%where $\displaystyle{\partial\rho^T}$ is the tangent component of $\partial \rho$;
%\item[(2)]
For any star-shaped $\mm_0$ the following evolution equation holds:
\begin{eqnarray*}
\frac{\partial Q(\mm_t)}{\partial t}&= & \frac{1}{n}Q(\mm_t)+|\mm_t|^{-1+\frac 1n}\int\left(\left(\frac{2n-1}{\siq}-\frac{1}{\coq}\right)\frac vH\right)d\mu\\
&&-|\mm_t|^{-1+\frac 1n}\int\left(\frac 1H\left(|A|^2-2(n+1)\right)\right)d\mu.
\end{eqnarray*}
%\end{itemize}
\end{Lemma}
\proof
%\begin{itemize}
%\item[(1)] 
Since $\hat{H}=(2n-1)\frac{\coh}{\sih}+\frac{\sih}{\coh}$ and $\frac{\partial\rho}{\partial t}=\frac vH$, it follows easly that
$$\displaystyle{\frac{\partial \hat{H}}{\partial t}=\frac{v}{H}\left(\frac{1}{\coq}-\frac{2n-1}{\siq}\right)}.$$
The thesis follows using this computation, the evolution of $H$ in Lemma \ref{eq_evoluz} and the fact that $$\frac{\Delta H}{H^2}-2\frac{|\nabla H|^2}{H^3}=-\Delta\left(\frac{1}{H}\right),$$ hence its integral vanishes.

\cvd

Now we want to show that if $Q$ decreases, it decreases very slowly.

\begin{Proposition}\label{evoluz_Q}
Let $\mm_t$ an $\sss^1$-invariant star-shaped hypersurface of $\mathbb{CH}^n$ evolving by inverse mean curvature flow. There is a positive constants $c$ which depends only on $n$ and $\mm_0$ such that
\begin{eqnarray*}
\frac{\partial Q(\mm_t)}{\partial t} & \geq & -ce^{-\frac{t}{n}}.%-ce^{\frac{t}{2n}}|\nabla H|.
\end{eqnarray*} 
\end{Proposition}
\proof By \eqref{2ff_s1}, \eqref{mc_s1} and \eqref{mc_s1_sigma} we can compute:

\begin{eqnarray}
\nonumber |A|^2-2(n+1) & = & h_i^kh_k^i-2(n+1)\\
\nonumber& = & \frac{\hat{\varphi}_{ij}\hat{\varphi}_{kh}\tilde{e}^{jk}\tilde{e}^{hi}}{v^2\siq}-\frac{2\coh}{v^2\siq}\hat{\varphi}_{ij}\tilde{e}^{ji}-2(n+1)\\
\nonumber&&(2n-1)\frac{\coq}{v^2\siq}+\frac{\siq}{v^2\coq}+\frac{2}{v^2}\\
\nonumber & = & \frac{\hat{\varphi}_{ij}\hat{\varphi}_{kh}\tilde{e}^{jk}\tilde{e}^{hi}}{v^2\siq}-2\frac{\coh}{v\sih}\left(H-\hat{H}+{\hat{H}}\left(1-\frac{1}{v}\right)\right)\\%\label{d03}
\nonumber& & +\frac{2n-1}{v^2\siq}-\frac{1}{v^2\coq}-2(n+1)\left(1-\frac{1}{v^2}\right).
\end{eqnarray}
Moreover, by Lemma \ref{hessiani}
$$
\hat{\varphi}_{ij}\hat{\varphi}_{kh}\tilde{e}^{jk}\tilde{e}^{hi}={\varphi}_{ij}{\varphi}_{kh}\tilde{\sigma}^{jk}\tilde{\sigma}^{hi}+2\siq|\nabla_{\sigma}\varphi|^2_{\sigma}.
$$
Then, by Lemma \ref{eq_evoluz_2} we have:
\begin{eqnarray*}
|\mm_t|^{1-\frac{1}{n}}\frac{\partial Q}{\partial t}&= &\int\left(\frac{1}{n}-\frac{2\coh}{vH\sih}\right)\left(H-\hat{H}\right)d\mu-\int\frac{\varphi_{ij}\varphi_{kh}\tilde{\sigma}^{jk}\tilde{\sigma}^{ih}}{v^2H\siq}d\mu\\ 
& & +\int\frac{1}{H}\left(\frac{2n-1}{\siq}-\frac{1}{\coq}\right)\left(v-\frac{1}{v^2}\right)d\mu\\
%&&+\int\frac{1}{H}\left(2(n+1)-\frac{2\hat{H}\coh}{(1+v)\sih}\right)\left(1-\frac{1}{v^2}\right)d\mu.
&&+\int\frac{|\nabla_{\sigma}\varphi|^2_{\sigma}}{v^2H}\left(2n-\frac{2}{1+v}\frac{\coh}{\sih}\hat{H}\right)d\mu.
\end{eqnarray*}

\nid By Proposition \ref{hess_phi_exp}, Proposition \ref{gradphi}, Proposition \ref{H_lambda_exp} and the fact that the radius grows like $\frac{t}{2n}$, we can estimate all the terms in the evolution of $Q$ in the following way

\begin{eqnarray*}
\left|\frac 1n-\frac {2\coh}{vH\sih}\right| & \leq & \left|\frac 1n-\frac 2H\right|+\frac{2}{H}\left|1-\frac {\coh}{\sih}\right|+\frac{2\coh}{H\sih}\left|1-\frac 1v\right|\\
& \leq&  ce^{-\frac tn};\\
\left|H-{\hat{H}}\right|& \leq & \left|H-2n\right|+\left|2n-\hat{H}\right|\leq c e^{-\frac tn};\\
\left|\frac{\varphi_{ij}\varphi_{kh}\tilde{\sigma}^{jk}\tilde{\sigma}^{ih}}{v^2H\siq}\right|&\leq& ce^{-\frac{2t}{n}};\\
1-\frac{1}{v^2}&\leq&  ce^{-\frac{t}{n}};
\end{eqnarray*}
hence:
\begin{eqnarray*}
\left(\frac{2n-1}{\siq}-\frac{1}{\coq}\right)\left(v-\frac{1}{v^2}\right)&\geq & -\frac{1}{\coq}\frac{v^2+v+1}{v^2(v+1)}|\nabla_{\sigma}\varphi|^2_{\sigma}\geq -ce^{-\frac{2t}{n}}.
\end{eqnarray*}
Finally, 
\begin{eqnarray*}
\left|2n-\frac{2}{1+v}\frac{\coh}{\sih}\hat{H}\right| & \leq & \left|2n-\hat{H}\right|+\hat{H}\left|\frac{\coh}{\sih}-1\right|+\hat{H}\frac{\coh}{\sih}\left|1-\frac{2}{1+v}\right|\leq ce^{-\frac{t}{n}}.
\end{eqnarray*}
\noindent Therefore 
\begin{eqnarray*}
\frac{\partial Q}{\partial t} & \geq & -c|\mm_t|^{-1+\frac{1}{n}}\int e^{-\frac{2t}{n}}d\mu\geq -c e^{-\frac{t}{n}}.
\end{eqnarray*}
\cvd

Now, following the strategy of \cite{HW} we can complete the proof of Theorem \ref{main}

\begin{Proposition}
There is an $\mm_0$ such that the rescaled induced metric $\tilde{g}_{\infty}$ does not have constant Webster curvature.
\end{Proposition}
\proof
Fix a positive constant $c_0$ and let $f:\sss^{2n-1}\rightarrow\mathbb{R}$ be an $\sss^1$-invariant function  such that
$$
\left(\int_{\sss^{2n-1}}e^{2nf}d\sigma\right)^{-1+\frac 1n}\int_{\sss^{2n-1}}e^{2nf}\left(e^{-f}\Delta_{\sigma}e^{-f}-n|\nabla_{\sigma} e^{-f}|^2\right) d\sigma\geq 4c_0.
$$
Consider the family of $\sss^1$-invariant star-shaped hyperfurfaces $\widetilde{\mm}^{\tau}$ defined by the radial function $\tilde{\rho}^{\tau}(z)=\tau+f(z)$. We can fix a $\tau$ big enough such that $\widetilde{\mm}^{\tau}$ is mean convex, and, by Proposition \ref{LQ}, $Q(\widetilde{\mm}^{\tau})\geq 2{c_0}$. Let $\mm_t^{\tau}$ be the evolution by inverse mean curvature flow of such $\widetilde{\mm}^{\tau}$. We want to estimate the evolution of $Q({\mm}_t^{\tau})$. The constant $c$ appearing in Proposition \ref{evoluz_Q} depends on $n$ and the initial datum, hence on $n$, $f$ and $\tau$. In view of \eqref{stime tau}, in our case it can be written as
$$
c=\tilde{c}e^{-2\tau},
$$
where $\tilde{c}$ depends only on $n$ and $f$. It follows that, up to increase even more the parameter $\tau$, Proposition \ref{evoluz_Q} ensures that 
$$
\lim_{t\rightarrow\infty} Q(\mm_t^{\tau})\geq {c_0}>0.
$$
The thesis follows from Proposition \ref{LQ}.
\cvd

We finish noting that we can find an $\sss^1$-invariant function $f$ such that 
$$
\left(\int_{\sss^{2n-1}}e^{2nf}d\sigma\right)^{-1+\frac 1n}\int_{\sss^{2n-1}}e^{2nf}\left(e^{-f}\Delta_{\sigma}e^{-f}-n|\nabla_{\sigma} e^{-f}|_{\sigma}^2\right) d\sigma.
$$
is large as desired. We will show that such an example exists in $\sss^3$, but an analogous construction holds also in higher dimension. 

\begin{Example}
Consider $\sss^3$ immersed in $\mathbb{C}^2$. Let $(z_1,z_2)$ be its complex coordinates and $\zeta=|z_2|^2-|z_1|^2$. Note that $\zeta$ is an $\sss^1$-invariant function on $\sss^3$. For every $k\in\mathbb{N}$, let us define
$$
f_k: (z_1,z_2)\in\sss^3\mapsto k\zeta\in\mathbb R.
$$
Some explicit computations show that
$$
e^{4f_k}\left(e^{-f_k}\Delta_{\sigma}e^{-f_k}-2|\nabla_{\sigma} e^{-f_k}|_{\sigma}^2\right)=4ke^{2f_k}\left(k\zeta^2+2\zeta-k\right).
$$
It follows that there exist a constant $\gamma$ independent on $k$ such that 
$$
\begin{array}{rcl}
Q_k&:=&\left(\int_{\sss^3}e^{4f_k}d\sigma\right)^{-\frac 12}\int_{\sss^3} e^{4f_k}\left(e^{-f_k}\Delta_{\sigma}e^{-f_k}-2|\nabla_{\sigma} e^{-f_k}|_{\sigma}^2\right)d \sigma\\
 &=& \gamma k\left(\int_{-1}^1\sqrt{1-\zeta^2}e^{4k\zeta}d\zeta\right)^{-\frac 12}\int_{-1}^1\sqrt{1-\zeta^2}\left(k\zeta^2+2\zeta-k\right)e^{2k\zeta}d\zeta\\
 & = & \gamma k \left(\frac{\pi}{4k}I_1(4k)\right)^{-\frac{1}{2}}\frac{\pi}{4k}I_2(2k),
\end{array}
$$
where $I_p(x)$ is the modified Bessel function of the first kind. As $x$ goes to infinity we have the asymptotic expansion:
$$
I_p(x)\sim \frac{e^x}{\sqrt{2\pi x}}.
$$
It follows that 
$$
Q_k\sim\gamma'k^{\frac 14}\ \text{ as }\  k\rightarrow\infty,
$$
for some constant $\gamma'$  independent on $k$.

\end{Example}

\bigskip

\noindent Giuseppe Pipoli, Institut Fourier, Universit\'e Grenoble Alpes, France.\\ E-mail: giuseppe.pipoli@univ-grenoble-alpes.fr

%Second email: pipolig@libero.it
%Please, use this second email to communicate with me.


\begin{thebibliography}{abc99xys}
\bibitem[Be]{Be}
\textsc{A. L. BESSE,} \emph{Einstein manifolds} Springer-Verlag, Berlin, Hidelberg, New York, 1987.

%\bibitem[AB]{AB} 
%\textsc{B. ANREWS, C. BAKER}, \emph{Mean curvature flow of pinched submanifolds to spheres}, J. Differential Geometry  \textbf{85} (2010), 357-395.

\bibitem[Di]{Di}
\textsc{Q. DING,} \emph{The inverse mean curvature flow in rotationally symmetric spaces}, Chin. Ann. Math. \textbf{32B(1)} (2011), 27 - 44.

\bibitem[DT]{DT}
\textsc{S. DRAGOMIR, G. TOMASSINI,} \emph{Differential geometry and analysis on CR manifolds}, Progress in math. vol 246, Birkh\"auser (2006).

\bibitem[Ge1]{Ge1}
\textsc{C. GERHARDT,} \emph{Flow of nonconvex hypersurfaces into spheres}, J. Differential Geometry, \textbf{32} (1990), 299 - 314.

\bibitem[Ge2]{Ge2}
\textsc{C. GERHARDT,} \emph{Curvature problems} Series in Geometry and Topology, vol. 39, International Press, Somerville, MA, (2006).

\bibitem[Ge3]{Ge3}
\textsc{C. GERHARDT,} \emph{Invererse mean curvature flows in hyperbolic space,} J. Differential Geometry \textbf{89}, (2011), 487 - 527.

\bibitem[HI]{HI}
\textsc{G. HUISKEN, T. ILMANEN,} \emph{The inverse mean curvature flow and the Riemannian Penrose inequality}, J. Differential Geometry, \textbf{59(3)} (2001), 353 - 437.

\bibitem[Hu]{Hu}
\textsc{G. HUISKEN,} \emph{Contracting convex hypersurfaces in Riemannian manifolds by their mean curvature}, Invent. Math. \textbf{84(3)} (1986),  463 - 480.

\bibitem[HW]{HW}
\textsc{P.K. HUNG, M.T. WANG,} \emph{Inverse mean curvature flow in the hyperbolic 3-space revisited}, Calculus of variation and partial differential equations \textbf{54(1)} (2015), 119 - 126.

\bibitem[JL]{JL}
\textsc{D. JERISON, J. LEE,} \emph{Extremals for the Sobolev inequality on the Heisenberg group and the CR Yamabe problem} J. AMS \textsc{1}(1), (1988), 1 - 13.


\bibitem[KS]{KS}
\textsc{N. KOIKE, Y. SAKAI,} \emph{The inverse mean curvature flow in rank one symmetric spaces of non-compact type}, Kyushu J. Math. \textbf{69} (2015), 259 - 284.

\bibitem[Kr]{Kr}
\textsc{NIKOLAI VLADIMIROVICH KRYLOV,} \emph{Nonlinear elliptic and parabolic equations of the secondorder,} Mathematics and its Applications (Soviet Series), 7, D. Reidel Publishing Co., Dordrecht, (1987).


%\bibitem[LW]{LW}
%\textsc{H.H. LI, Y. WEI}, \emph{On inverse mean curvature flow in Schwarzschild space and Kottler space}, Calc. Var. (2017), 56 - 62.
\bibitem[Ne]{Ne}
\textsc{ANDR\'E NEVES,} \emph{Insufficient convergence of inverse mean curvature flow on asymptotically hyperbolic manifolds}, J. Differential Geometry \textbf{84} (2016), 191 - 229.

\bibitem[O]{O}
\textsc{B. O'NEILL}, \emph{The fundamental equations of a submersion} Michigan Math. J. \textbf{13} (1966) 459 - 469.

\bibitem[Pa]{Pa}
\textsc{J.R. PARKER,} \emph{Notes on complex hyperbolic geometry}, http://maths.dur.ac.uk/~dma0jrp/img/NCHG.pdf (2003).

\bibitem[Pi1]{Pi}
\textsc{G. PIPOLI,} \emph{Mean curvature flow and Riemannian submersions}, Geom. Dedicata \textbf{184(1)} (2016), 67 - 81.

\bibitem[Pi2]{Pi2}
\textsc{G. PIPOLI,} \emph{Inverse mean curvature flow in quaternionic hyperbolic space}, to appear on Atti Accad. Naz. Lincei Rend. Lincei Mat. Appl.

\bibitem[Sc]{Sc}
\textsc{J. SCHEUER,} \emph{The inverse mean curvature flow in warped cylinders of non-positive radial curvature}, Advances in Mathematics \textbf{306} (2017), 1130 - 1163.

\bibitem[Ur]{Ur}
\textsc{J.I.E. URBAS,} \emph{On the expansion of starshaped hypersurfaces by symmetric functions of their principal curvatures,} Math. Z., \textbf{205}, (1990), 355 - 372.

\bibitem[Zh]{Zh}
\textsc{H. ZHOU,} \emph{Inverse mean curvature flows in warped product manifolds} to appear on The Journal of Geometric Analysis.


%\bibitem[XW]{XW}
%\textsc{X. Wang,} \emph{On a remarkable formula of Jerison and Lee in CR geometry} Math. Res. Lett., vol. 22, \textbf{01} (2015), 279 - 299.
\end{thebibliography}
\end{document}